\documentclass[12pt]{article}
\usepackage{amsfonts,amssymb,subeqnarray,eqsection,indent}
\usepackage{graphicx}

\begin{document}

\title{\vspace*{-2.5cm} Zero-free Regions \\
       for Multivariate Tutte Polynomials \\
       (alias Potts-model Partition Functions) \\
       of Graphs and Matroids}

\author{
  {\small Bill Jackson}                                    \\[-2mm]
  {\small\it School of Mathematical Sciences}  \\[-2mm]
  {\small\it Queen Mary, University of London} \\[-2mm]
  {\small\it Mile End Road} \\[-2mm]
  {\small\it London E1 4NS, England}                         \\[-2mm]
  {\small\tt B.JACKSON@QMUL.AC.UK}                        \\[5mm]
  {\small Alan D.~Sokal}\thanks{Also at Department of Mathematics,
           University College London, London WC1E 6BT, England.}  \\[-2mm]
  {\small\it Department of Physics}       \\[-2mm]
  {\small\it New York University}         \\[-2mm]
  {\small\it 4 Washington Place}          \\[-2mm]
  {\small\it New York, NY 10003 USA}      \\[-2mm]
  {\small\tt SOKAL@NYU.EDU}               \\[-2mm]
  {\protect\makebox[5in]{\quad}}  
}

\date{June 19, 2008 \\[5mm]
      \emph{Dedicated to the memory of Bill Tutte (1917--2002)}}

\maketitle
\thispagestyle{empty}   

\vspace*{-5mm}

\begin{abstract}
The chromatic polynomial $P_G(q)$ of a loopless graph $G$
is known to be nonzero (with explicitly known sign)
on the intervals $(-\infty,0)$, $(0,1)$ and $(1,32/27]$.
Analogous theorems hold for the flow polynomial of bridgeless graphs
and for the characteristic polynomial of loopless matroids.
Here we exhibit all these results as special cases of more general theorems
on real zero-free regions of the multivariate Tutte polynomial
$Z_G(q,{\bf v})$.
The proofs are quite simple, and employ deletion-contraction
together with parallel and series reduction.
In particular, they shed light on the origin of the curious number 32/27.
\end{abstract}


\bigskip
\noindent
{\bf Mathematics Subject Classification (MSC) codes:}
05C15 (Primary); 05A20, 05B35, 05C99, 05E99, 82B20 (Secondary).

\clearpage

\newtheorem{defin}{Definition}[section]
\newtheorem{definition}{Definition}[section]
\newtheorem{prop}[defin]{Proposition}
\newtheorem{proposition}[defin]{Proposition}
\newtheorem{lem}[defin]{Lemma}
\newtheorem{lemma}[defin]{Lemma}
\newtheorem{guess}[defin]{Conjecture}
\newtheorem{ques}[defin]{Question}
\newtheorem{question}[defin]{Question}
\newtheorem{prob}[defin]{Problem}
\newtheorem{problem}[defin]{Problem}
\newtheorem{thm}[defin]{Theorem}
\newtheorem{theorem}[defin]{Theorem}
\newtheorem{cor}[defin]{Corollary}
\newtheorem{corollary}[defin]{Corollary}
\newtheorem{conj}[defin]{Conjecture}
\newtheorem{conjecture}[defin]{Conjecture}

\newtheorem{pro}{Problem}
\newtheorem{clm}{Claim}
\newtheorem{con}{Conjecture}

%
%
\newcounter{example}[section]
\newenvironment{example}%
{\refstepcounter{example}
 \bigskip\par\noindent{\bf Example \thesection.\arabic{example}.}\quad
}%
{\quad $\Box$}
\def\bexam{\begin{example}}
\def\eexam{\end{example}}

\renewcommand{\theenumi}{\alph{enumi}}
\renewcommand{\labelenumi}{(\theenumi)}
\def\prf{\par\noindent{\bf Proof.\enspace}\rm}
\def\rmk{\par\medskip\noindent{\bf Remark.\enspace}\rm}

\newcommand{\be}{\begin{equation}}
\newcommand{\ee}{\end{equation}}
\newcommand{\<}{\langle}
\renewcommand{\>}{\rangle}
\newcommand{\widebar}{\overline}
\def\reff#1{(\protect\ref{#1})}
\def\spose#1{\hbox to 0pt{#1\hss}}
\def\ltapprox{\mathrel{\spose{\lower 3pt\hbox{$\mathchar"218$}}
 \raise 2.0pt\hbox{$\mathchar"13C$}}}
\def\gtapprox{\mathrel{\spose{\lower 3pt\hbox{$\mathchar"218$}}
 \raise 2.0pt\hbox{$\mathchar"13E$}}}
\def\textprime{${}^\prime$}
\def\proof{\par\medskip\noindent{\sc Proof.\ }}
\def\sketchofproof{\par\medskip\noindent{\sc Sketch of Proof.\ }}
\newcommand{\qed}{\quad $\Box$ \medskip \medskip}
\def\proofof#1{\bigskip\noindent{\sc Proof of #1.\ }}
\def\half{ {1 \over 2} }
\def\third{ {1 \over 3} }
\def\twothird{ {2 \over 3} }
\def\smfrac#1#2{\textstyle{#1\over #2}}
\def\smhalf{ \smfrac{1}{2} }
\newcommand{\real}{\mathop{\rm Re}\nolimits}
\renewcommand{\Re}{\mathop{\rm Re}\nolimits}
\newcommand{\imag}{\mathop{\rm Im}\nolimits}
\renewcommand{\Im}{\mathop{\rm Im}\nolimits}
\newcommand{\sgn}{\mathop{\rm sgn}\nolimits}
\def\hboxscript#1{ {\hbox{\scriptsize\em #1}} }

\newcommand{\restrict}{\upharpoonright}
\renewcommand{\emptyset}{\varnothing}

\newcommand{\scra}{{\mathcal{A}}}
\newcommand{\scrb}{{\mathcal{B}}}
\newcommand{\scrc}{{\mathcal{C}}}
\newcommand{\scrf}{{\mathcal{F}}}
\newcommand{\scrg}{{\mathcal{G}}}
\newcommand{\scrh}{{\mathcal{H}}}
\newcommand{\scrl}{{\mathcal{L}}}
\newcommand{\scrm}{{\mathcal{M}}}
\newcommand{\scro}{{\mathcal{O}}}
\newcommand{\scrp}{{\mathcal{P}}}
\newcommand{\scrr}{{\mathcal{R}}}
\newcommand{\scrs}{{\mathcal{S}}}
\newcommand{\scrt}{{\mathcal{T}}}
\newcommand{\scrv}{{\mathcal{V}}}
\newcommand{\scrw}{{\mathcal{W}}}
\newcommand{\scrz}{{\mathcal{Z}}}

\newcommand{\w}{{\bf w}}
\newcommand{\wtilde}{{\widetilde{\bf w}}}
\newcommand{\what}{{\widehat{\bf w}}}
\newcommand{\z}{{\bf z}}
\newcommand{\Rtilde}{{\widetilde{\bf R}}}
\newcommand{\Rhat}{{\widehat{\bf R}}}
\newcommand{\K}{{\bf K}}
\newcommand{\p}{{\bf p}}
\renewcommand{\k}{{\bf k}}
\newcommand{\n}{{\bf n}}
\renewcommand{\a}{{\bf a}}
\renewcommand{\b}{{\bf b}}
\renewcommand{\r}{{\bf r}}
\newcommand{\smalln}{\mbox{\scriptsize\bf n}}
\newcommand{\blambda}{\mbox{\boldmath $\lambda$}}
\newcommand{\smallblambda}{\mbox{\scriptsize\boldmath $\lambda$}}
\newcommand{\btheta}{\mbox{\boldmath $\theta$}}
\newcommand{\balpha}{\mbox{\boldmath $\alpha$}}
\newcommand{\bdelta}{\mbox{\boldmath $\delta$}}
\newcommand{\smallbalpha}{\mbox{\scriptsize\boldmath $\alpha$}}
\def\twotilde#1{{\widetilde{\widetilde{#1}}}}
\def\twohat#1{{\widehat{\widehat{#1}}}}

\newcommand{\C}{{\mathbb C}}
\newcommand{\Z}{{\mathbb Z}}
\newcommand{\N}{{\mathbb N}}
\newcommand{\R}{{\mathbb R}}
\newcommand{\Cbar}{{\overline{\C}}}

\newcommand{\Zhat}{\widehat{Z}}
\newcommand{\Ztilde}{\widetilde{Z}}
\newcommand{\Ptilde}{\widetilde{P}}
\newcommand{\Ctilde}{\widetilde{C}}
\newcommand{\series}{{\,\bowtie_q\,}}
\newcommand{\seriesq}{{\,\bowtie_q\,}}
\newcommand{\seriesnoq}{{\,\bowtie\,}}
\renewcommand{\parallel}{\Vert}

\newcommand{\bk}{ {\bf k} }
\newcommand{\bm}{ {\bf m} }
\newcommand{\bp}{ {\bf p} }
\newcommand{\bu}{ {\bf u} }
\newcommand{\bv}{ {\bf v} }
\newcommand{\bw}{ {\bf w} }
\newcommand{\bx}{ {\bf x} }
\newcommand{\bzero}{ {\bf 0} }
\def\qvbf{{q, \bv}}
\def\veff{v_{{\rm eff}}}

\def\zgxconny{Z_G^{(x \leftrightarrow y)}}
\def\zgxnoconny{Z_G^{(x \not\leftrightarrow y)}}

%
%
\font\fourrm  = cmr5 
\def\sspose#1{\hbox to 0pt{#1\hss}}
\def\gtalmost{\mathrel{\sspose{\lower 0.75pt\hbox{\kern-1.5pt\fourrm (\quad)}}
  \raise 2.0pt\hbox{$\ge$}}}
\def\ltalmost{\mathrel{\sspose{\lower 0.75pt\hbox{\kern-1.5pt\fourrm (\quad)}}
  \raise 2.0pt\hbox{$\le$}}}


\newenvironment{sarray}{
	  \textfont0=\scriptfont0
	  \scriptfont0=\scriptscriptfont0
	  \textfont1=\scriptfont1
	  \scriptfont1=\scriptscriptfont1
	  \textfont2=\scriptfont2
	  \scriptfont2=\scriptscriptfont2
	  \textfont3=\scriptfont3
	  \scriptfont3=\scriptscriptfont3
	\renewcommand{\arraystretch}{0.7}
	\begin{array}{l}}{\end{array}}

\newenvironment{scarray}{
	  \textfont0=\scriptfont0
	  \scriptfont0=\scriptscriptfont0
	  \textfont1=\scriptfont1
	  \scriptfont1=\scriptscriptfont1
	  \textfont2=\scriptfont2
	  \scriptfont2=\scriptscriptfont2
	  \textfont3=\scriptfont3
	  \scriptfont3=\scriptscriptfont3
	\renewcommand{\arraystretch}{0.7}
	\begin{array}{c}}{\end{array}}



\section{Introduction}

It is known (see e.g.\ \cite{Jackson_03})
that the chromatic polynomial $P_G(q)$ of a loopless graph $G$ satisfies:
\begin{theorem}
  \label{thm1.1}
Let $G$ be a loopless graph that has
$n$ vertices, $c$ components, and $b$ nontrivial blocks.
[We call a block ``trivial'' if it has only one vertex,
 and ``nontrivial'' otherwise.]
Then:
\begin{itemize}
   \item[(a)]  $P_G(q)$ is nonzero with sign $(-1)^n$ for $q \in (-\infty,0)$.
   \item[(b)]  $P_G(q)$ has a zero of multiplicity $c$ at $q=0$.
   \item[(c)]  $P_G(q)$ is nonzero with sign $(-1)^{n+c}$ for $q \in (0,1)$.
   \item[(d)]  $P_G(q)$ has a zero of multiplicity $b$ at $q=1$.
   \item[(e)]  $P_G(q)$ is nonzero with sign $(-1)^{n+c+b}$ for
       $q \in (1,{32 \over 27}]$.
\end{itemize}
\end{theorem}
Analogous theorems are also known for the flow polynomial of bridgeless graphs
and, more generally, for the characteristic polynomial of loopless matroids.

\begin{figure}[t]
  \vspace*{-5mm}
  \centering
  \includegraphics[width=5.5in]{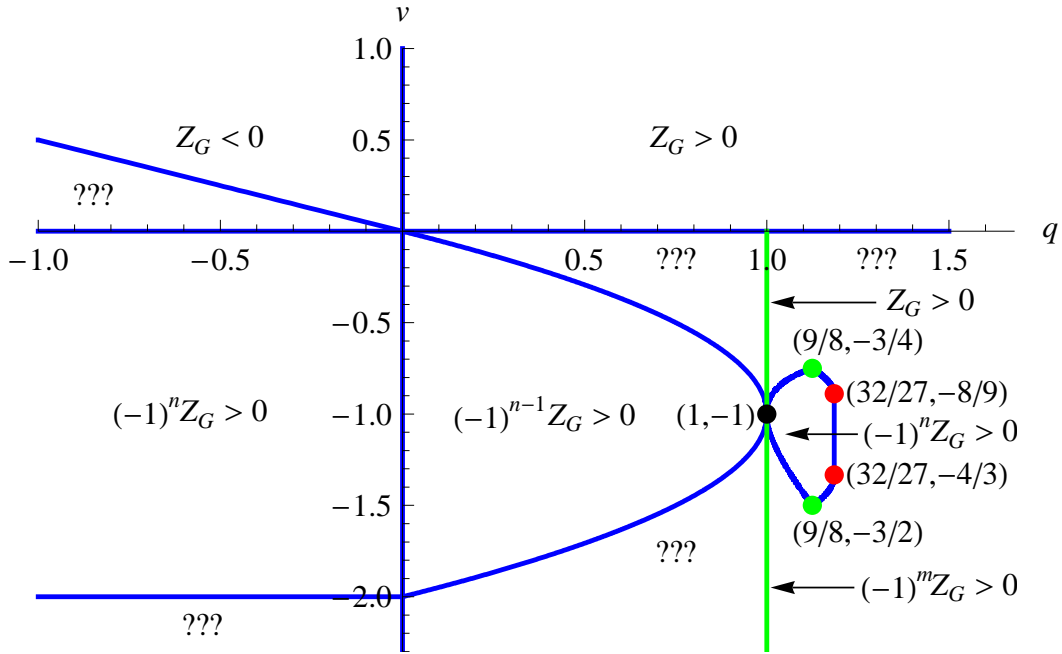}
  \caption{
      Sign of $Z_G(q,{\bf v})$ for (loopless) 2-connected graphs $G$.
      The cases $q>0$, $v \ge 0$ and $q=1$ are trivial
      (see Section~\ref{sec.trivial}).
      The remaining cases $q<0$, $0 < q < 1$ and $1 < q \le 32/27$
      are proven in Sections~\ref{sec.q<0}, \ref{sec.0<q<1}
      and \ref{sec.q>1}, respectively.
  }
\label{fig.introplot}
\end{figure}

All the foregoing polynomials are special cases of the
multivariate Tutte polynomial $Z_G(q,{\bf v})$
--- also known as the Potts-model partition function in statistical
mechanics --- or its generalization to matroids
(see \cite{Sokal_bcc2005} for a recent survey).
Here ${\bf v} = \{v_e\}_{e \in E}$ are real or complex edge weights,
and one recovers the chromatic (resp.\ flow) polynomial
if one sets $v_e = -1$ (resp.\ $v_e = -q$) for all edges $e$.
The purpose of this paper is to exhibit all the results
of types (a), (b), (c) and (e) as special cases of theorems on
real zero-free regions of the multivariate Tutte polynomial.
Our results are illustrated in Figure~\ref{fig.introplot}.

%
%

One message of the present paper (see also \cite{Sokal_bcc2005})
is that there is considerable advantage
in studying the multivariate polynomial $Z_G(q, \bv)$,
even if one is ultimately interested in a particular
two-variable or one-variable specialization.
For instance, $Z_G(q, \bv)$ is {\em multiaffine}\/
in the variables $\bv$ (i.e., of degree 1 in each $v_e$ separately);
and often a multiaffine polynomial in many variables is easier to handle
than a general polynomial in a single variable
(e.g., it may permit simple proofs by induction on the number of variables).
Furthermore, many natural operations on graphs,
such as the reduction of edges in series or parallel,
lead out of the class of ``all $v_e$ equal''.
For these reasons, the multivariate extension
of a single-variable result is sometimes much {\em easier}\/
to prove than its single-variable specialization.
Examples of the advantage obtained by considering general $\{v_e\}$ are:
\begin{itemize}
   \item[(a)] a simple proof of the Brown--Hickman theorem
       on chromatic roots of large subdivisions
       (compare \cite[Appendix A]{Sokal_chromatic_roots}
         with \cite{Brown-Hickman_99b});
   \item[(b)] a very simple proof of the (multivariate) Brown--Colbourn
       conjecture on the zeros of reliability polynomials,
       in the special case of series-parallel graphs
       (compare \cite[Remark 3 in Section 4.1]{Sokal_chromatic_bounds}
        with \cite{Wagner_00}); and
   \item[(c)] a disproof of the Brown--Colbourn conjecture
       for general graphs \cite{Royle-Sokal}.
       (Both the univariate and multivariate conjectures are
       false for general graphs;  but a counterexample to the
       univariate conjecture would have been very difficult to find
       by direct search.  Rather, one first shows that
       the complete graph $K_4$ is a counterexample to the
       multivariate conjecture;  one then uses the formulae for
       parallel connection of edges to find a 16-edge counterexample
       to the univariate conjecture.)
\end{itemize}
In this paper we shall give further examples of the utility of
considering general $\{v_e\}$;  in particular, we shall elucidate
the origin of the curious number 32/27 in Theorem~\ref{thm1.1}(e).
A further advantage of the $Z_G(q,\bv)$ formalism
is that it shows clearly the distinct roles
played by the variables $q$ and $\{v_e\}$:
namely, $q$ is a global parameter
while the edge weights $\{v_e\}$ are variables that can be mapped.

A second message of this paper is that it is sometimes advantageous
to ``think matroidal'', even when the ultimate goal is to study graphs.
Indeed, as Oxley \cite{Oxley_01} has eloquently shown,
graph theorems can often be improved by rethinking them
in matroidal terms --- that is, by eliminating reference to concepts
that have no matroidal analogue
(e.g.\ vertices and their degrees, connected components, \ldots)
and replacing them by matroidal concepts
(e.g.\ rank, circuits, cocircuits, \ldots).
Another advantage of working with matroids is that
every matroid has a dual, while only {\em planar}\/ graphs
have duals with reasonable algebraic properties.
The matroidal philosophy is particularly pertinent in the present case,
because the multivariate Tutte polynomial can be defined naturally
for matroids (Section~\ref{sec2.1})
and even in the graphical case it ``sees'' only the underlying
matroidal structure
(that is, two graphs with the same cycle matroid
have the same multivariate Tutte polynomial,
modulo trivial factors of $q$).
For this reason, we believe that matroids are the ``natural'' category
for studying the multivariate Tutte polynomial.

The plan of this paper is as follows:
In Section~\ref{sec2} we review the definition of the
multivariate Tutte polynomial for graphs and matroids,
along with some of its elementary properties.
In Section~\ref{sec.trivial} we briefly discuss the trivial cases
$q > 0$, $v \ge 0$ and $q=1$.
In Sections~\ref{sec.q<0} and \ref{sec.0<q<1}
we study the intervals $q < 0$ and $0 < q < 1$, respectively.
In Section~\ref{sec.abstract} we prove an abstract result
that will be important in what follows.
In Section~\ref{subsec.0<q<1.block} we strengthen the results for $0 < q < 1$
by considering the block structure of $G$.
In Section~\ref{sec.diamond} we collect some properties of the
``diamond map'', which plays a fundamental role in our analysis.
In Section~\ref{sec.q>1} we study the interval $1 < q \le 32/27$.
Finally, in Section~\ref{sec.further} we state some conjectured
extensions of our results.
%

Since some readers of this paper may be unfamiliar with matroids,
we shall ordinarily state and prove each theorem first for graphs
and only afterwards for matroids, even though logically speaking
the latter contains the former.  In most cases the matroidal proofs
will be nearly direct translations of the graphical proofs
into matroidal language;  we shall therefore usually be brief
in discussing the matroidal proofs, drawing attention only to
any non-obvious points.

\section{The multivariate Tutte polynomial}  \label{sec2}

\subsection{Definition for graphs and matroids}  \label{sec2.1}

Let $G=(V,E)$ be a finite undirected graph
with vertex set $V$ and edge set $E$;
in this paper all graphs are allowed to have
loops and multiple edges unless explicitly stated otherwise.
The {\em multivariate Tutte polynomial}\/ of $G$ is, by definition,
the polynomial
\begin{equation}
   Z_G(q, \bv)   \;=\;
   \sum_{A \subseteq E}  q^{k(A)}  \prod_{e \in A}  v_e
   \;,
  \label{eq1.1.bcc2005}
\end{equation}
where $q$ and $\bv = \{v_e\}_{e \in E}$ are commuting indeterminates,
and $k(A)$ denotes the number of connected components in the subgraph $(V,A)$.
[It is sometimes convenient to consider instead
\begin{equation}
   \Ztilde_G(q, \bv)   \;\equiv\;
   q^{-|V|} Z_G(q, \bv)  \;=\;
   \sum_{A \subseteq E}  q^{k(A) - |V|}  \prod_{e \in A}  v_e
   \;,
  \label{eq.def.Ztilde}
\end{equation}
which is a polynomial in $q^{-1}$ and $\{v_e\}$.]
{}From a combinatorial point of view,
$Z_G$ is simply the multivariate generating polynomial
that enumerates the spanning subgraphs of $G$ according to
their precise edge content (with weight $v_e$ for the edge~$e$)
and their number of connected components (with weight $q$ for each component).
As we shall see, $Z_G$ encodes a vast amount of
combinatorial information about the graph $G$,
and contains many other well-known graph polynomials as special cases.
In this paper we shall take an analytic point of view,
and treat $q$ and $\{v_e\}$ as real variables.\footnote{
   The study of $Z_G(q, \bv)$ when $q$ and $\{v_e\}$
   are treated as {\em complex}\/ variables is also of
   great interest to both mathematicians and statistical physicists:
   see e.g.\ \cite{Sokal_chromatic_bounds,Sokal_chromatic_roots,%
Shrock_BCC99,Salas-Sokal_transfer1,Jackson_03,Royle-Sokal,Sokal_bcc2005}.
}

If we set all the edge weights $v_e$ equal to the same value $v$,
we obtain a two-variable polynomial $Z_G(q,v)$ that is
equivalent to the standard Tutte polynomial $T_G(x,y)$
after a simple change of variables
[see \reff{eq.Tutte.2} below].

All of these considerations can be extended from graphs
to matroids.\footnote{
   See Oxley \cite{Oxley_92} for an excellent introduction to
   matroid theory.
}
Let $M$ be a matroid with ground set $E$
and rank function $r_M \colon\, 2^{E} \to \N$.
We then define the multivariate Tutte polynomial
\be
   \Ztilde_M(\qvbf)  \;=\;
   \sum_{A \subseteq E}  q^{-r_M(A)}  \prod_{e \in A}  v_e
   \;,
 \label{eq.defZ.matroid}
\ee
which is a polynomial in $q^{-1}$ and $\{v_e\}$.
This extends the graph definition \reff{eq.def.Ztilde}
in the sense that if $G$ is a graph and $M(G)$ is its cycle matroid, then
\begin{equation}
   \Ztilde_{M(G)}(q, \bv)   \;=\;
   \Ztilde_G(q, \bv)
\end{equation}
[because $r_{M(G)}(A) = |V| - k(A)$].
Since a matroid is completely determined by its rank function,
$\Ztilde_M$ is simply an algebraic encoding of {\em all}\/
the information about the matroid $M$.
Moreover, our earlier statement that $Z_G$ encodes ``a vast amount''
of information about the graph $G$ can now be made more precise:
$Z_G$ encodes the number of vertices $|V|$
together with all the information about $G$
that is contained in its cycle matroid $M(G)$ [and no other information].
In particular, if $G$ is loopless and 3-connected,
then it is uniquely determined
(within the class of loopless graphs without isolated vertices)
by its cycle matroid $M(G)$
and hence by $\Ztilde_G$ (or equivalently by $Z_G$);
this is a special case of Whitney's 2-isomorphism theorem
\cite[Theorem~5.3.1 and Lemma~5.3.2]{Oxley_92}.

Let us also remark \cite{Sokal_bcc2005} that the
multivariate Tutte polynomial of a matroid $M$ is related
to that of its dual matroid $M^*$ by the formula
\begin{subeqnarray}
   \Ztilde_{M^*}(q, \bv)  & = &
   q^{-r_{M^*}(E)} \left( \prod\limits_{e \in E} v_e \right)
                                                          \Ztilde_M(q, q/\bv)
       \\[2mm]
   & = &
   q^{r_M(E)} \left( \prod\limits_{e \in E} {v_e \over q} \right)
                                                          \Ztilde_M(q, q/\bv)
   \;.
 \label{eq.duality.matroid}
\end{subeqnarray}
(Here $q/\bv=\{q/v_e\}_{e\in E}$,   $r_M(E)$ is the rank of $M$,
$r_{M^*}(E)$ is the rank of $M^*$, and we have 
$r_M(E)+r_{M^*}(E)=|E|$.)
In~brief, duality takes $v_e \mapsto q/v_e$
(and inserts some prefactors).
Indeed, the duality formula \reff{eq.duality.matroid}
is an easy consequence of the definition \reff{eq.defZ.matroid}
together with the formula for the rank function of a dual:
\begin{equation}
   r_{M^*}(A)  \;=\;  |A| \,+\, r_M(E \setminus A) \,-\, r_M(E)
   \;.
\end{equation}
It goes without saying that the duality formula \reff{eq.duality.matroid}
can be specialized from matroids to {\em planar}\/ graphs.
One of the advantages of working with matroids
is that we can think about duality even for non-planar graphs.

It is convenient to introduce explicitly the coefficients of $Z_G(\qvbf)$
as a polynomial in $q$:
\be
   Z_G(\qvbf)   \;=\;   \sum\limits_{k=1}^n C_G^{[k]}({\bf v}) \, q^k
 \label{eq.potts.coeffs.1}
\ee
where $n=|V|$ and
\be
   C_G^{[k]}({\bf v})   \;=\;
   \sum\limits_{\begin{scarray}
                   A \subseteq E \\
                   k(A) = k
                \end{scarray}}
      \prod_{e \in A}  v_e
   \;.
 \label{eq.potts.coeffs.defCGk}
\ee
Likewise, let us introduce the coefficients of $\Ztilde_M(\qvbf)$
as a polynomial in $q^{-1}$:
\be
   \Ztilde_M(\qvbf)   \;=\;
   \sum\limits_{r=0}^{r(M)} \Ctilde_M^{[r]}({\bf v}) \, q^{-r}
 \label{eq.matroid.coeffs.1}
\ee
where
\be
   \Ctilde_M^{[r]}({\bf v})   \;=\;
   \sum\limits_{\begin{scarray}
                   A \subseteq E \\
                   r_M(A) = r
                \end{scarray}}
      \prod_{e \in A}  v_e
   \;.
 \label{eq.matroid.coeffs.defCGk}
\ee

\subsection{Coloring interpretation for graphs}

Let $G=(V,E)$ be a finite graph and let $q$ be a positive integer.
A {\em proper $q$-coloring}\/ of $G$
is a map $\sigma \colon\, V \to \{ 1,2,\ldots,q \}$
such that $\sigma(i) \neq \sigma(j)$ for all pairs of adjacent vertices $i,j$.
It is not hard to show (see below)
that for each graph $G$ there exists a polynomial $P_G(q)$
with integer coefficients such that, for each $q \in \Z_+$,
the number of proper $q$-colorings of $G$ is precisely $P_G(q)$.
This (obviously unique) polynomial $P_G(q)$
is called the {\em chromatic polynomial}\/ of $G$.\footnote{
   See \cite{Read_68,Read_88} for excellent reviews on chromatic polynomials,
   and \cite{Chia_97} for an extensive bibliography through 1995.
   A review of results and conjectures
   concerning the zeros of chromatic polynomials
   can be found in \cite{Jackson_03}.
}

A more general polynomial can be obtained as follows:
Assign to each edge $e \in E$ a real or complex weight $v_e$,
and write ${\bf v} = \{v_e\}_{e \in E}$ for the collection
of these weights.
Then the {\em $q$-state Potts-model partition function}\/
for the graph $G$ is defined by
\begin{equation}
   Z_G^{\rm Potts}(\qvbf)   \;=\;
   \sum_{ \sigma \colon\, V \to \{ 1,2,\ldots,q \} }
   \; \prod_{e \in E}  \,
      \biggl[ 1 + v_e \delta(\sigma_{x_1(e)}, \sigma_{x_2(e)}) \biggr]
   \;.
 \label{def.ZPotts}
\end{equation}
Here the sum runs over all maps $\sigma\colon\, V \to \{ 1,2,\ldots,q \}$,
and we sometimes write $\sigma_x$ as a synonym for $\sigma(x)$;
the $\delta$ is the Kronecker delta
\begin{equation}
   \delta(a,b)   \;=\;   \cases{1  & if $a=b$ \cr
                                \noalign{\vskip 2pt}
                                0  & if $a \neq b$ \cr
                               }
\end{equation}
and $x_1(e), x_2(e) \in V$ are the two endpoints of the edge $e$
(in arbitrary order).
In particular, if we take $v_e = -1$ for all $e$,
then a coloring $\sigma$ gets weight 1 or 0
according as it is proper or improper,
so that $Z_G^{\rm Potts}(q,-1)$ counts the proper $q$-colorings.

In statistical physics, the formula \reff{def.ZPotts} arises as follows:
In the Potts model \cite{Potts_52,Wu_82,Wu_84},
an ``atom'' (or ``spin'') at the site $x \in V$ can exist in any one of
$q$ different states.
A {\em configuration}\/ is a map $\sigma\colon\, V \to \{1,\ldots,q\}$.
The {\em energy}\/ of a configuration is the sum, over all edges $e \in E$,
of $0$ if the spin values at the two endpoints of that edge are unequal
and $-J_e$ if they are equal.
The {\em Boltzmann weight}\/ of a configuration is then $e^{-\beta H}$,
where $H$ is the energy of the configuration
and $\beta \ge 0$ is the inverse temperature.
The {\em partition function}\/ is the sum, over all configurations,
of their Boltzmann weights.
Clearly this is just a rephrasing of \reff{def.ZPotts},
with $v_e = e^{\beta J_e} - 1$.
A parameter value $J_e$ (or $v_e$)
is called {\em ferromagnetic}\/ if $J_e \ge 0$ ($v_e \ge 0$),
as it is then favored for adjacent spins to take the same value;
{\em antiferromagnetic}\/ if $-\infty \le J_e \le 0$ ($-1 \le v_e \le 0$),
as it is then favored for adjacent spins to take different values;
and {\em unphysical}\/ if $v_e \notin [-1,\infty)$,
as the weights are then no longer nonnegative.
The chromatic polynomial ($v_e=-1$) thus corresponds
to the zero-temperature ($\beta \to\infty$) limit
of the antiferromagnetic ($J_e < 0$) Potts model.
The main idea of the present paper is that many results
for chromatic polynomials extend to part or all of the
antiferromagnetic regime
(and indeed into part of the unphysical regime as well).

It is far from obvious that $Z_G^{\rm Potts}(\qvbf)$,
which is defined separately for each positive integer $q$,
is in fact the restriction to $q \in \Z_+$
of a {\em polynomial}\/ in $q$.
But this is in fact the case, and indeed we have:

\begin{theorem}[Fortuin--Kasteleyn \protect\cite{Kasteleyn_69,Fortuin_72}
    representation of the Potts model]
   \label{thm.FK}
For integer $q \ge 1$,
\begin{equation}
   Z_G^{\rm Potts}(\qvbf) \;=\;  Z_G(q, \bv)  \;.
 \label{eq.FK.identity}
\end{equation}
That is, the Potts-model partition function
is simply the specialization of the multivariate Tutte polynomial
to $q \in \Z_+$.
\end{theorem}

\proof
In \reff{def.ZPotts}, expand out the product over $e \in E$,
and let $A \subseteq E$ be the set of edges for which the term
$v_e \delta(\sigma_{x_1(e)}, \sigma_{x_2(e)})$ is taken.
Now perform the sum over configurations $\{ \sigma_x \}_{x \in V}$:
in each component of the subgraph $(V,A)$
the color $\sigma_x$ must be constant,
and there are no other constraints.
Therefore,
\begin{equation}
   Z^{\rm Potts}_G(\qvbf) \;=\;
   \sum_{ A \subseteq E }  q^{k(A)}  \prod_{e \in A}  v_e
   \;,
  \label{eq1.1bis}
\end{equation}
as was to be proved.
\qed

The subgraph expansion \reff{eq1.1bis} was discovered by
Birkhoff \cite{Birkhoff_12} and Whitney \cite{Whitney_32a}
for the special case $v_e = -1$ (see also Tutte \cite{Tutte_47,Tutte_54});
in its general form it is due to
Fortuin and Kasteleyn \cite{Kasteleyn_69,Fortuin_72}
(see also \cite{Edwards-Sokal}).

Special cases of the multivariate Tutte polynomial $Z_G(\qvbf)$
include the chromatic polynomial ($v=-1$) and the flow polynomial ($v=-q$),
and more generally the standard two-variable Tutte polynomial.
Indeed, we have
\begin{eqnarray}
   P_G(q)   & = &   Z_G(q,-1)   \\[2mm]
   F_G(q)   & = &   q^{-|V|} (-1)^{|E|} Z_G(q,-q)   \\[2mm]
   T_G(x,y) & = &   (x-1)^{-k(E)} \, (y-1)^{-|V|} \,
                         Z_G \Bigl( (x-1)(y-1), \, y-1 \Bigr)   \;.
 \label{eq.Tutte.2}
\end{eqnarray}
Several other evaluations of the multivariate Tutte polynomial
are discussed in \cite{Sokal_bcc2005}.

\subsection{Elementary identities}

We now wish to prove some elementary identities
for the multivariate Tutte polynomial.
There are two alternative approaches to proving such identities:
one is to prove the identity directly for real or complex $q$
(or considering $q$ as an algebraic indeterminate),
using the subgraph expansion \reff{eq1.1.bcc2005}
or its generalization \reff{eq.defZ.matroid} to matroids;
the other is to prove the identity first for {\em positive integer}\/ $q$,
using the coloring representation \reff{def.ZPotts}/\reff{eq.FK.identity},
and then to extend it to general $q$ by arguing that two polynomials
(or rational functions) that coincide at infinitely many points must be equal.
The latter approach is perhaps less elegant,
but it is often simpler or more intuitive.
However, only the former approach extends to arbitrary matroids.

One way to guess (albeit not to prove) an identity for matroids
is to prove it first for graphs,
and then translate it from $Z_G$ to $\Ztilde_G = q^{-|V|} Z_G$;
usually the latter identity carries over verbatim to matroids,
{\em mutatis mutandis}\/.

In this paper we shall use four principal tools:
factorization over blocks,
the deletion-contraction identity,
the parallel-reduction identity,
and the series-reduction identity.

\bigskip

{\bf Factorization.}\/
If $G$ is the disjoint union of $G_1$ and $G_2$, then trivially
\begin{equation}
   Z_G(q,\bv)  \;=\;  Z_{G_1}(q,\bv) \, Z_{G_2}(q,\bv)
   \;.
 \label{eq.components}
\end{equation}
That is, $Z_G$ ``factorizes over components''.

A slightly less trivial situation arises when $G$ consists of
subgraphs $G_1$ and $G_2$ joined at a single cut vertex $x$;
in this case
\begin{equation}
   Z_G(q,\bv)  \;=\;  {Z_{G_1}(q,\bv) \, Z_{G_2}(q,\bv)  \over q}
   \;.
 \label{eq.blocks}
\end{equation}
This is easily seen from the subgraph expansion in the variant form
\begin{equation}
   Z_G(q, \bv)   \;=\;
   q^{|V|} \, \sum_{A \subseteq E}  q^{\gamma(A)} \prod_{e \in A} {v_e \over q}
   \;,
 \label{eq1.1cycles}
\end{equation}
where
\begin{equation}
   \gamma(A)   \;=\;   k(A) \,-\, |V| \,+\, |A|
 \label{eq.cyclomatic}
\end{equation}
is the cyclomatic number
(i.e., number of linearly independent cycles) of the graph $(V,A)$.
It is also easily seen from the coloring representation
\reff{def.ZPotts}/\reff{eq.FK.identity}
by first fixing the color $\sigma_x$ at the cut vertex
and then summing over it;
from this viewpoint, \reff{eq.blocks} reflects
the $S_q$ permutation symmetry of the $q$-state Potts model.\footnote{
   More precisely, it reflects the symmetry of the spin model
   under a global transformation $\sigma_y \mapsto g\sigma_y$
   (simultaneously for all $y \in V$)
   that acts {\em transitively}\/ on each single-spin space.
}
We summarize \reff{eq.blocks} by saying that
$Z_G$ ``factorizes over blocks'' modulo a factor $q$.

The identities \reff{eq.components} and \reff{eq.blocks}
can be written in a unified form, by using $\Ztilde_G = q^{-|V|} Z_G$:
in both cases we have
\begin{equation}
   \Ztilde_G(q,\bv)  \;=\;  \Ztilde_{G_1}(q,\bv) \, \Ztilde_{G_2}(q,\bv)
   \;.
 \label{eq.factorize.ZtildeG}
\end{equation}
This, in turn, is a special case of the following obvious fact:
if a matroid $M$ is the direct sum of matroids $M_1$ and $M_2$, then
\begin{equation}
   \Ztilde_M(q,\bv)  \;=\;  \Ztilde_{M_1}(q,\bv) \, \Ztilde_{M_2}(q,\bv)
   \;.
\end{equation}

\bigskip

{\bf Deletion-contraction identity.}\/
If $e \in E$, let $G \setminus e$
denote the graph obtained from $G$ by deleting the edge $e$,
and let $G/e$ denote the graph obtained from $G \setminus e$
by contracting the two endpoints of $e$ into a single vertex
(please note that we retain in $G/e$ any loops or multiple edges
 that may be formed as a result of the contraction).
Then, for any $e \in E$, we have the identity
\be
   Z_G(\qvbf)   \;=\;
      Z_{G \setminus e}(q, {\bf v}_{\neq e})   \,+\,
      v_e Z_{G/e}(q, {\bf v}_{\neq e})
   \;.
 \label{eq.delcon}
\ee
This is easily seen either from the coloring representation
\reff{def.ZPotts}/\reff{eq.FK.identity}
or the subgraph expansion \reff{eq1.1.bcc2005}.
Please note that the deletion-contraction identity \reff{eq.delcon}
takes the same form regardless of whether $e$ is a normal edge,
a loop, or a bridge
(in contrast to the situation for the usual Tutte polynomial $T_G$).
Of course, if $e$ is a loop, then $G/e = G \setminus e$,
so we can also write
\be
   Z_G  \;=\;  (1+v_e) Z_{G \setminus e}  \;=\;  (1+v_e) Z_{G/e}
   \qquad\hbox{if $e$ is a loop}\;.
 \label{eq.delcon.loop}
\ee
Similarly, if $e$ is a bridge, then $G \setminus e$ is the disjoint union
of two subgraphs $G_1$ and $G_2$ while $G/e$ is obtained by joining
$G_1$ and $G_2$ at a cut vertex,
so that $Z_{G/e} = Z_{G \setminus e}/q$ and hence
\be
   Z_G  \;=\;  (1+v_e/q) Z_{G \setminus e}  \;=\;  (q+v_e) Z_{G/e}
   \qquad\hbox{if $e$ is a bridge}\;.
 \label{eq.delcon.bridge}
\ee

The deletion-contraction identity applies also to
the coefficients $C_G^{[k]}$ of the multivariate Tutte polynomial:
\be
   C_G^{[k]}({\bf v})   \;=\;
      C_{G \setminus e}^{[k]}({\bf v}_{\neq e})   \,+\,
      v_e C_{G/e}^{[k]}({\bf v}_{\neq e})
   \;.
 \label{eq.potts.coeffs.delcon_CGk}
\ee
This follows either by examining the definition \reff{eq.potts.coeffs.defCGk}
or by observing that the deletion-contraction identity \reff{eq.delcon}
for $Z_G$ does not mix powers of $q$.

In terms of $\Ztilde_G = q^{-|V|} Z_G$,
the deletion-contraction identity takes the form
\begin{subeqnarray}
   \Ztilde_G  & = &  \Ztilde_{G \setminus e} \,+\, {v_e \over q} \Ztilde_{G/e}
                        \qquad\hbox{if $e$ is not a loop}   \\[3mm]
   \Ztilde_G  & = &  \Ztilde_{G \setminus e} \,+\, v_e \Ztilde_{G/e}
                                                       \nonumber \\
              & = &  (1+v_e) \Ztilde_{G \setminus e}   \nonumber \\
              & = &  (1+v_e) \Ztilde_{G/e}
                        \qquad\quad\;\hbox{if $e$ is a loop}
 \label{eq.delcon.Ztilde}
\end{subeqnarray}
as easily follows from \reff{eq.delcon}
together with the counting of vertices in $G \setminus e$ and $G/e$.

Not surprisingly, the deletion-contraction formula for matroids
is identical in form to \reff{eq.delcon.Ztilde}:
\begin{subeqnarray}
   \Ztilde_M  & = &  \Ztilde_{M \setminus e} \,+\, {v_e \over q} \Ztilde_{M/e}
                        \qquad\hbox{if $e$ is not a loop}   \\[3mm]
   \Ztilde_M  & = &  \Ztilde_{M \setminus e} \,+\, v_e \Ztilde_{M/e}
                                                       \nonumber \\
              & = &  (1+v_e) \Ztilde_{M \setminus e}   \nonumber \\
              & = &  (1+v_e) \Ztilde_{M/e}
                        \qquad\quad\;\hbox{if $e$ is a loop}
 \label{eq.delcon.matroid}
\end{subeqnarray}
This easily follows from the formulae for the rank function
of a deletion or contraction:  if $A \subseteq E \setminus e$, then
\begin{subeqnarray}
   r_{M \setminus e}(A)   & = &   r_M(A)   \\[3mm]
   r_{M/e}(A)             & = &
        \cases{
            r_M(A \cup e) - 1  & if $e$ is not a loop \cr
            r_M(A \cup e)      & if $e$ is a loop     \cr
        }
  \slabel{eq.rank.contraction}
\end{subeqnarray}

\bigskip

{\bf Parallel-reduction identity.}\/
If $G$ contains edges $e_1,e_2$
connecting the same pair of vertices $x,y$, they can be replaced,
without changing the value of $Z$,
by a single edge $e=xy$ with weight
\be
   v_e  \;=\;  (1+v_{e_1})(1+v_{e_2}) \,-\, 1
        \;=\;  v_{e_1} + v_{e_2} + v_{e_1} v_{e_2}
   \;.
 \label{eq.parallel1}
\ee
This is easily seen either from the coloring representation
\reff{def.ZPotts}/\reff{eq.FK.identity}
or the subgraph expansion \reff{eq1.1.bcc2005}.
More formally, we can write
\be
   Z_G(q, {\bf v}_{\neq e_1,e_2}, v_{e_1}, v_{e_2})  \;=\;
   Z_{G \setminus e_2}(q, {\bf v}_{\neq e_1,e_2},
                          v_{e_1} + v_{e_2} + v_{e_1} v_{e_2})
   \;.
 \label{eq.parallel2}
\ee
The parallel-reduction rule
$(v_1,v_2) \mapsto v_{\rm eff}$ with $1 + v_{\rm eff} = (1+v_1)(1+v_2)$
can be remembered by the mnemonic ``$1+v$ multiplies''.
We write $v_1 \parallel v_2 \equiv (1+v_1)(1+v_2) - 1$;
and if $\scrv_1,\scrv_2 \subseteq \R$ (or $\C$)
we write $\scrv_1 \parallel \scrv_2 \equiv \{ v_1 \parallel v_2 \colon\,
   v_1 \in \scrv_1,\, v_2 \in \scrv_2 \}$.

The parallel-reduction rule applies also to the $C_G^{[k]}$:
\be
   C_G^{[k]}({\bf v}_{\neq e_1,e_2}, v_{e_1}, v_{e_2})  \;=\;
   C_{G \setminus e_2}^{[k]}({\bf v}_{\neq e_1,e_2},
                          v_{e_1} + v_{e_2} + v_{e_1} v_{e_2})
   \;.
 \label{eq.parallel_CGk}
\ee
This follows either from the definition \reff{eq.potts.coeffs.defCGk}
or by observing that the parallel-reduction rule
for $Z_G$ does not mix powers of $q$.

A virtually identical formula holds for matroids:
if $e_1$ and $e_2$ are parallel elements in a matroid $M$
(i.e., form a two-element circuit), then
\begin{equation}
   \Ztilde_M(q, \bv_{\neq e_1,e_2}, v_{e_1}, v_{e_2})  \;=\;
   \Ztilde_{M \setminus e_2}(q, \bv_{\neq e_1,e_2},
                          v_{e_1} + v_{e_2} + v_{e_1} v_{e_2})
   \;.
 \label{eq.parallel3}
\end{equation}
The formula \reff{eq.parallel3} also holds trivially if
$e_1$ and $e_2$ are both loops.

\bigskip

{\bf Series-reduction identity.}\/
We say that edges $e_1, e_2 \in E$
are {\em in series (in the narrow sense)}\/
if there exist vertices $x,y,z \in V$ with $x \neq y$ and $y \neq z$
such that $e_1$ connects $x$ and $y$,
$e_2$ connects $y$ and $z$, and $y$ has degree 2 in $G$.
In this case the pair of edges $e_1,e_2$ can be replaced,
without changing the value of $Z$,
by a single edge $e=xz$ with weight
\be
   v_e  \;=\;  {v_{e_1} v_{e_2}  \over  q + v_{e_1} + v_{e_2}}
 \label{eq.series1}
\ee
provided that we then multiply $Z$ by the prefactor $q + v_{e_1} + v_{e_2}$.
More formally, we can write
\be
   Z_G(q, {\bf v}_{\neq e_1,e_2}, v_{e_1}, v_{e_2})  \;=\;
   (q + v_{e_1} + v_{e_2}) \,
   Z_{G / e_2}(q, {\bf v}_{\neq e_1,e_2},
                          v_{e_1} v_{e_2} / (q + v_{e_1} + v_{e_2}))
   \;.
 \label{eq.series2}
\ee
This identity can be derived from the coloring representation
\reff{def.ZPotts}/\reff{eq.FK.identity}
by noting that
\begin{subeqnarray}
   \sum_{\sigma_y = 1}^q
   [1 + v_{e_1} \delta(\sigma_x,\sigma_y)]
   [1 + v_{e_2} \delta(\sigma_y,\sigma_z)]
   & = &
   q + v_{e_1} + v_{e_2} + v_{e_1} v_{e_2} \delta(\sigma_x,\sigma_z)  \\
   & = &
   (q + v_{e_1} + v_{e_2})
   \left[1 \,+\, {v_{e_1} v_{e_2}  \over  q + v_{e_1} + v_{e_2}}
                 \delta(\sigma_x,\sigma_z)
   \right]
   \;.
     \nonumber \\
\end{subeqnarray}
Alternatively, it can be derived from the subgraph expansion
\reff{eq1.1.bcc2005}
by considering the four possibilities for the edges $e_1$ and $e_2$
to be occupied or empty and analyzing the number of connected components
thereby created.
The series-reduction rule
$(v_1,v_2) \mapsto v_{\rm eff} \equiv v_1 v_2/(q+v_1+v_2)$
can be remembered by the mnemonic ``$1+q/v$ multiplies'':
namely,
\be
   1 + {q \over v_{\rm eff}}  \;=\;
   \left( 1 + {q \over v_1} \right)
   \left( 1 + {q \over v_2} \right)
   \;.
  \label{eq.series2a}
\ee
We write $v_1 \series v_2 \equiv v_1 v_2/(q+v_1+v_2)$;
and if $\scrv_1,\scrv_2 \subseteq \R$ (or $\C$)
we write $\scrv_1 \series \scrv_2 \equiv \{ v_1 \series v_2 \colon\,
   v_1 \in \scrv_1,\, v_2 \in \scrv_2 \}$.

Consider now the more general situation in which
$\{e_1,e_2\}$ is a two-edge cut of $G$
(not necessarily the cut associated with a degree-2 vertex $y$);
we then say that $e_1,e_2$ are {\em in series (in the wide sense)}\/.
It turns out that the identity \reff{eq.series2} still holds.
To see this, let us prove the generalization of this identity to matroids.
Let $e_1$ and $e_2$ be series elements in a matroid $M$,
i.e., suppose that $\{e_1,e_2\}$ is a cocircuit.
Then, for any $A \subseteq E \setminus \{e_1,e_2\}$, we have
\begin{equation}
   r_M(A \cup e_1)  \;=\; r_M(A \cup e_2)  \;=\;  r_M(A) + 1
 \label{eq.rank.cocircuit}
\end{equation}
(since the complement of a cocircuit is a hyperplane).
A short calculation using \reff{eq.rank.contraction} with $e=e_2$ then yields
\begin{equation}
   \Ztilde_M(q, \bv_{\neq e_1,e_2}, v_{e_1}, v_{e_2})  \;=\;
   {q + v_{e_1} + v_{e_2} \over  q}
   \,
   \Ztilde_{M / e_2}(q, \bv_{\neq e_1,e_2},
                          v_{e_1} v_{e_2} / (q + v_{e_1} + v_{e_2}))
   \;.
 \label{eq.series3}
\end{equation}
%
%
The formula \reff{eq.series3} also holds trivially if
$e_1$ and $e_2$ are both coloops.

Please note that duality $v \mapsto q/v$ interchanges
the parallel-reduction rule (``$1+v$ multiplies'')
with the series-reduction rule (``$1+q/v$ multiplies'').
This is no accident,
since we now see that parallel-reduction and series-reduction
(in the wide sense)
are indeed duals of each other:
$\{e_1,e_2\}$ is a circuit (resp.\ cocircuit) in $M$
if and only if it is a cocircuit (resp.\ circuit) in the dual matroid $M^*$.

\section{Two trivial cases}   \label{sec.trivial}

Let us begin by disposing of two cases in which we can trivially
control the sign of $Z_G(q,\bv)$.

\begin{proposition}
   \label{prop.trivial1}
Let $G$ be a graph with $n$ vertices, and suppose that $v_e \ge 0$
for all $e$.  Then $Z_G(q,\bv) \ge q^n > 0$ for $q > 0$;
and more generally, for $0 \le \ell \le n$ we have
\be
   {d^\ell \over dq^\ell} \, Z_G(q,\bv)
   \;\ge\; n^{\underline \ell} q^{n-\ell}  \;>\; 0  \qquad\hbox{for } q > 0
\ee
[here $n^{\underline \ell} = n(n-1) \cdots (n-\ell+1)$].
\end{proposition}

\proof
In the definition \reff{eq1.1.bcc2005},
the term $A = \emptyset$ contributes $q^n$,
and the remaining terms contribute a polynomial in $q$
with nonnegative coefficients.
\qed

A similar result holds for matroids, but since $\Ztilde_M(\qvbf)$
involves inverse powers of $q$ [cf.\ \reff{eq.defZ.matroid}],
we can no longer control the derivatives with respect to $q$:

\begin{proposition}
   \label{prop.trivial1.matroid}
Let $M$ be a matroid, and suppose that $v_e \ge 0$ for all $e$.
Then $\Ztilde_M(q,\bv) \ge 1 > 0$ for $q > 0$.
\end{proposition}

The other trivial case is $q=1$, because
$Z_G(1,\bv) = \prod\limits_{e \in E} (1+v_e)$
and more generally
$\Ztilde_M(\qvbf) = \prod\limits_{e \in E} (1+v_e)$.
It follows that:

\begin{proposition}
   \label{prop.trivial2}
Let $G$ be a graph with $m$ edges.
\begin{itemize}
   \item[(a)]  If $v_e > -1$ for all $e$, then $Z_G(1,\bv) > 0$.
   \item[(b)]  If $v_e < -1$ for all $e$, then $(-1)^m Z_G(1,\bv) > 0$.
\end{itemize}
\end{proposition}

\noindent
The corresponding result of course holds for a matroid
on a ground set with $m$ elements.

\section{The interval \boldmath $q\in (-\infty,0)$}   \label{sec.q<0}

It is well known that the coefficients of the chromatic polynomial
alternate in sign,  and  the leading coefficient is 1
whenever the graph is loopless.
These facts immediately imply that
$P_G(q)$ is nonzero with sign $(-1)^n$ for $q < 0$.
The following theorem generalizes this result
to the multivariate Tutte polynomial $Z_G(\qvbf)$:

\begin{theorem}
   \label{thm.q<0}
Let $G$ be a graph with $n$ vertices and $c$ components.
Suppose that
\begin{itemize}
   \item[(i)] $v_e \ge -1$ for every loop $e$; and
   \item[(ii)] $-2 \le v_e \le 0$ for every non-loop edge $e$.
\end{itemize}
Then:
\begin{itemize}
   \item[(a)] $C_G^{[k]} \equiv 0$ for $1 \le k < c$,
      $(-1)^{n-k} C_G^{[k]}({\bf v}) \ge 0$ for $c \le k \le n$,
      and \\
$C_G^{[n]} = \prod\limits_{{\rm loops}\, e} (1+v_e)$.
   \item[(b)]  $(-1)^n Z_G(\qvbf) \ge 0$ for all $q < 0$,
      with strict inequality if and only if every loop has $v_e > -1$.
\end{itemize}
Furthermore, if
\begin{itemize}
   \item[(i${}'$)] $v_e > -1$ for every loop $e$; and
   \item[(ii${}'$)] $-2 < v_e < 0$ for every non-loop edge $e$,
\end{itemize}
then:
\begin{itemize}
   \item[(a${}'$)] $(-1)^{n-k} C_G^{[k]}({\bf v}) > 0$ for $c \le k \le n$.
   \item[(b${}'$)] The root of $Z_G(\qvbf)$ at $q=0$ has multiplicity $c$.
\end{itemize}
\end{theorem}

Please note that Theorem~\ref{thm1.1}(a,b)
are the special cases of Theorem~\ref{thm.q<0}(b,b${}'$)
in which the graph is loopless and $v_e = -1$ for all edges $e$.
We now see that these results can be generalized to
$v_e \in [-2,0]$ and $v_e \in (-2,0)$, respectively.
In particular, Theorem~\ref{thm1.1}(a,b) extends to the
whole antiferromagnetic regime $[-1,0]$
as well as to part of the unphysical regime $(-\infty,-1]$.

Let us also remark that since $q < 0$, conclusion (b)
can equivalently be written as $\Ztilde_G(\qvbf) \ge 0$.
This way of writing Theorem~\ref{thm.q<0}
also suggests that the correct generalization to matroids
will be $\Ztilde_M(\qvbf) \ge 0$:
see Theorem~\ref{thm.q<0.matroid} below.

\smallskip

\proofof{Theorem~\ref{thm.q<0}}
It suffices to prove (a), since (b) is then an immediate corollary.
Since each loop $e$ simply contributes an overall factor $1 + v_e$,
which has the right sign by hypothesis,
we can assume henceforth that $G$ is loopless.
By \reff{eq.potts.coeffs.defCGk}, $C_G^{[k]} \equiv 0$ for $1 \le k < c$ and,
since $G$ is loopless, $C_G^{[n]} \equiv 1$.
The proof of the sign inequalities for $C_G^{[k]}$
is by induction on the number of edges in $G$.
If $G$ has no edges, then $c=n$ and (a) holds.
Now suppose that $G$ has $m$ edges,
and assume that the result holds for all graphs having fewer than $m$ edges.
We now consider three cases:

(i) If $G$ is a forest, then
we have $k(A) = n - |A|$ for every $A \subseteq E$, so that
\be
   C_G^{[k]}({\bf v})   \;=\;
   \sum\limits_{\begin{scarray}
                   A \subseteq E \\
                   |A| = n-k
                \end{scarray}}
      \prod_{e \in A}  v_e
   \;.
\ee
Since $v_e \le 0$ for all $e$, statement (a) holds.

We can henceforth suppose that $G$ has at least one circuit.

(ii)  If $G$ has somewhere a pair $e_1,e_2$ of parallel edges,
pick some such pair and apply the parallel-reduction formula
\reff{eq.parallel_CGk} to it.
Since $(v_1,v_2) \mapsto (1+v_1)(1+v_2) - 1$
maps the interval $[-2,0]$ into itself,
we can apply the inductive hypothesis to $G \setminus e_2$;
and the result has the right sign, since $G \setminus e_2$
has the same number of vertices and components as $G$ does.

(iii) If $G$ has no pair of parallel edges,
pick any edge $e$ which belongs to a circuit of $G$
and apply the deletion-contraction identity
\reff{eq.potts.coeffs.delcon_CGk}.
By the inductive hypothesis,
the first term has the right sign, since $G \setminus e$
has the same number of vertices and components as $G$ does.
Since $G$ has no parallel edges,
$G/e$ is loopless, so we can apply the inductive
hypothesis to it as well;
moreover, since $e$ is not a loop, $G/e$ has one vertex fewer than $G$
and the same number of components as $G$;
therefore, since $v_e \le 0$,
the second term has the right sign as well.

This proves (a);  and the same argument, with minor modifications,
proves (a${}'$) under the hypotheses (i${}'$) and (ii${}'$).
Statement (b${}'$) then follows from $C_G^{[k]}=0$ for $k<c$
and $C_G^{[k]} \neq 0$ for $k=c$.
\qed

The interval $[-2,0]$ is best possible, as is shown by the
following example:

\bexam
   \label{sec.q<0.exam1}
Let $G=K_2^{(m)}$ (a single pair of vertices connected by
$m$ parallel edges).
Then $Z_G(q,v) = q[q + (1+v)^m - 1]$,
so that for any $q < 0$ there are roots $v$
tending to $-2$ (from below) and to $0$ (from above)
as $m \to \infty$ (the former only for $m$ even).
\eexam


\bigskip
\bigskip
\bigskip


We can also prove some inequalities on the partial derivatives of
$C_G^{[k]}(\bv)$ with respect to individual weights $v_e$,
provided that we make a slightly stronger hypothesis on the interval
in which the weights $v_e$ lie.
If $e_1,\ldots,e_\ell$ are edges in $G$
and $e \in E \setminus \{e_1,\ldots,e_\ell\}$,
let us say that {\em $e$ is spanned by $\{e_1,\ldots,e_\ell\}$}\/
if there exists a subset of $\{e_1,\ldots,e_\ell\}$
which together with $e$ forms a circuit
[or equivalently, if $\{e_1,\ldots,e_\ell\} \cup \{e\}$
 has a cyclomatic number larger than that of $\{e_1,\ldots,e_\ell\}$;
 or equivalently, if the rank of $\{e_1,\ldots,e_\ell\} \cup \{e\}$
 in the cycle matroid $M(G)$ is equal to that of $\{e_1,\ldots,e_\ell\}$].
Note that a loop is spanned by any set of edges (even the empty set).

\begin{corollary}
     \label{cor_alternating_2}
Let $G=(V,E)$ be a graph with $n$ vertices,
and let $\ell\ge 0$ and $e_1,\ldots,e_\ell \in E$.
Suppose that
\begin{itemize}
   \item[(i)] $-2 \le v_e \le 0$ for all
     $e \in E \setminus \{e_1,\ldots,e_\ell\}$
     that are not spanned by $\{e_1,\ldots,e_\ell\}$; and
   \item[(ii)] $v_e \ge -1$ for all
     $e \in E \setminus \{e_1,\ldots,e_\ell\}$
     that are spanned by $\{e_1,\ldots,e_\ell\}$.
\end{itemize}
Then:
\begin{itemize}
      \item[(a)]  If $e_1,\ldots,e_\ell$ are not all distinct,
         we have
\be\label{2.26}
       {\partial^\ell C_G^{[k]}(\bv) \over
        \partial v_{e_1} \cdots \partial v_{e_\ell}}
       \;=\;  0   \;.
\ee
      \item[(b)]  If $e_1,\ldots,e_\ell$ are all distinct
and form
a subgraph with cyclomatic number $\gamma$, we have
\be\label{2.27}
       (-1)^{n-k+\ell+\gamma} \,
       {\partial^\ell C_G^{[k]}(\bv) \over
        \partial v_{e_1} \cdots \partial v_{e_\ell}}
       \;\ge\;  0   \;.
\ee
\end{itemize}
\end{corollary}

Please note that Theorem~\ref{thm.q<0}(a) is simply the
special case $\ell = 0$ of Corollary~\ref{cor_alternating_2}.

\proof
Note first that the deletion-contraction identity 
\reff{eq.potts.coeffs.delcon_CGk}
implies that
\be
      {\partial C_G^{[k]}(\bv) \over \partial v_e}  \;=\;
      C_{G/e}^{[k]}(\bv_{\neq e})
\ee
and more generally
\be
      {\partial^\ell C_G^{[k]}(\bv) \over
       \partial v_{e_1} \cdots \partial v_{e_\ell}}
      \;=\;
      C_{G/\{e_1,\ldots,e_\ell\}}^{[k]}(\bv_{\neq e_1,\ldots,e_\ell})
 \label{eq.derivCG.multiple}
\ee
for any set $e_1,\ldots,e_\ell$ of distinct edges.
[If the edges $e_1,\ldots,e_\ell$ are not distinct,
 then $\partial^\ell C_G^{[k]}(\bv) /
       \partial v_{e_1} \cdots \partial v_{e_\ell}  \equiv 0$.]
Now, the graph $G/e$ has $|V|$ vertices if $e$ is a loop,
and $|V|-1$ vertices if $e$ is not a loop.
Iterating this fact, we see that if $e_1,\ldots,e_\ell$ are distinct edges
that form a subgraph with cyclomatic number $\gamma$,
then $G/\{e_1,\ldots,e_\ell\}$ has $|V|-\ell+\gamma$ vertices.
Furthermore, if $e  \in E \setminus \{e_1,\ldots,e_\ell\}$,
then $e$ is a loop in $G/\{e_1,\ldots,e_\ell\}$
if and only if $e$ is spanned by $\{e_1,\ldots,e_\ell\}$.
The hypothesis $v_e \ge -1$ for such edges
is precisely what is needed to apply Theorem~\ref{thm.q<0}.
Applying \reff{eq.derivCG.multiple}
together with Theorem~\ref{thm.q<0}(a) and the foregoing observations
yields \reff{2.27}.
\qed

Corollary~\ref{cor_alternating_2} was proven a few years ago by
Scott and Sokal \cite[Proposition~2.7]{Scott-Sokal},
under the slightly stronger hypothesis that $-1 \le v_e \le 0$
for all $e \in E$.
However, their proof used a fairly sophisticated device,
namely the partitionability identity.
It is nice to know that a completely elementary proof can be given,
and that the conditions on $v_e$ can be slightly weakened.

\bigskip

The proof of Theorem~\ref{thm.q<0} extends immediately to matroids, yielding:

\begin{theorem}
   \label{thm.q<0.matroid}
Let $M$ be a matroid of rank $r(M)$ on the ground set $E$.
Suppose that
\begin{itemize}
   \item[(i)] $v_e \ge -1$ for every loop $e$; and
   \item[(ii)] $-2 \le v_e \le 0$ for every non-loop edge $e$.
\end{itemize}
Then:
\begin{itemize}
   \item[(a)]  $(-1)^{r} \Ctilde_M^{[r]}({\bf v}) \ge 0$ for $0 \le r \le r(M)$,
      with $\Ctilde_M^{[0]} = \prod\limits_{{\rm loops}\, e} (1+v_e)$.
   \item[(b)]  $\Ztilde_M(\qvbf) \ge 0$ for all $q < 0$,
      with strict inequality if and only if every loop has $v_e > -1$.
\end{itemize}
Furthermore, if
\begin{itemize}
   \item[(i${}'$)] $v_e > -1$ for every loop $e$; and
   \item[(ii${}'$)] $-2 < v_e < 0$ for every non-loop edge $e$,
\end{itemize}
then $(-1)^{r} \Ctilde_M^{[r]}({\bf v}) > 0$ for $0 \le r \le r(M)$.
\end{theorem}

Applying Theorem~\ref{thm.q<0.matroid} to the cographic matroid $M^*(G)$,
we obtain a result dual to Theorem~\ref{thm.q<0}:

\begin{theorem}[Dong \protect\cite{Dong_08}]
   \label{thm.q<0,dual}
Let $G$ be a graph with $c$ components, and fix $q < 0$.
Suppose that
\begin{itemize}
   \item[(i)] $v_e \ge -q$ for every bridge $e$; and
   \item[(ii)] $v_e \ge -q/2$ for every non-bridge edge $e$.
\end{itemize}
Then $(-1)^c  Z_G(\qvbf) \ge 0$,
with strict inequality if and only if every bridge has $v_e > -q$.
\end{theorem}

\noindent
This result was suggested to us by Feng-Ming Dong \cite{Dong_08},
who proved it by a direct argument that is essentially the dual
of the proof of Theorem~\ref{thm.q<0}.

Similarly, the proof of Corollary~\ref{cor_alternating_2}
extends immediately to matroids
(using the same definition of ``spanned'',
 which is after all the matroidal one).
We obtain:

\begin{corollary}
     \label{cor_alternating_2.matroid}
Let $M$ be a matroid of rank $r(M)$ on the ground set $E$,
and let $\ell\ge 0$ and $e_1,\ldots,e_\ell \in E$.
Suppose that
\begin{itemize}
   \item[(i)] $-2 \le v_e \le 0$ for all
     $e \in E \setminus \{e_1,\ldots,e_\ell\}$
     that are not spanned by $\{e_1,\ldots,e_\ell\}$; and
   \item[(ii)] $v_e \ge -1$ for all
     $e \in E \setminus \{e_1,\ldots,e_\ell\}$
     that are spanned by $\{e_1,\ldots,e_\ell\}$.
\end{itemize}
Then:
\begin{itemize}
      \item[(a)]  If $e_1,\ldots,e_\ell$ are not all distinct,
         we have
\be\label{2.26.matroid}
       {\partial^\ell \Ctilde_M^{[r]}(\bv) \over
        \partial v_{e_1} \cdots \partial v_{e_\ell}}
       \;=\;  0   \;.
\ee
      \item[(b)]  If $e_1,\ldots,e_\ell$ are all distinct
and form a set with rank $\rho$ in the matroid $M$, we have
\be\label{2.27.matroid}
       (-1)^{r+\rho} \,
       {\partial^\ell \Ctilde_M^{[r]}(\bv) \over
        \partial v_{e_1} \cdots \partial v_{e_\ell}}
       \;\ge\;  0   \;.
\ee
\end{itemize}
\end{corollary}


\section{The interval \boldmath $q\in (0,1)$}   \label{sec.0<q<1}

In this section we discuss the conditions under which
the sign of $Z_G(\qvbf)$ can be controlled
when $0 < q < 1$ and the edge weights $v_e$ lie in
a suitable subinterval of $(-2,0)$.
We
prove a basic result
valid for arbitrary graphs~$G$ (Theorem~\ref{thm.0<q<1.improved}).
Later, in  Section~\ref{sec4.2},
we will
prove a sequence of refinements that make successively stronger hypotheses
on the minimum number of edges in each block of $G$,
and obtain correspondingly wider intervals for the edge weights $v_e$.


Let $G$ be a loopless graph with $n$ vertices and $c$ components;
then Theorem~\ref{thm1.1}(c) states that
$P_G(q)$ is nonzero with sign $(-1)^{n+c}$ for $0 < q < 1$.
The following theorem generalizes this result
to the multivariate Tutte polynomial $Z_G(\qvbf)$:

\begin{theorem}
   \label{thm.0<q<1.improved}
Let $G$ be a graph with $n$ vertices and $c$ components,
and let $q \in (0,1)$.
Suppose that:
\begin{itemize}
   \item[(i)]  $v_e > -1$ for every loop $e$;
   \item[(ii)]  $v_e < -q$ for every bridge $e$; and
   \item[(iii)]  $-1 - \sqrt{1-q} < v_e < -1 + \sqrt{1-q}$
       for every normal (i.e., non-loop non-bridge) edge $e$.
\end{itemize}
Then $(-1)^{n+c} Z_G(\qvbf) > 0$.
\end{theorem}

\begin{corollary}
   \label{cor.0<q<1.improved}
Let $G$ be a graph with $n$ vertices and $c$ components,
and let $q \in (0,1)$.
Then $(-1)^{n+c} Z_G(\qvbf) > 0$ in each of the following three cases:
\begin{itemize}
   \item[(a)]  $G$ is loopless, with $-1 - \sqrt{1-q} < v_e < -q$
       for all $e \in E$.
   \item[(b)]  $G$ is bridgeless, with $-1 < v_e < -1 + \sqrt{1-q}$
       for all $e \in E$.
   \item[(c)]  $G$ is loopless and bridgeless,
       with $-1 - \sqrt{1-q} < v_e < -1 + \sqrt{1-q}$ for all $e \in E$.
\end{itemize}
\end{corollary}

Theorem~\ref{thm1.1}(c) is the special case of
Corollary~\ref{cor.0<q<1.improved}(a)
in which $v_e = -1$ for all edges $e$.
We now see that this result can be extended to
$v_e \in (-1-\sqrt{1-q},-q)$.
Please note that this interval approaches $(-2,0)$ as $q \downarrow 0$,
and degenerates to the empty set $(-1,-1)$ as $q \uparrow 1$.

It is worth remarking that the hypotheses of
Theorem~\ref{thm.0<q<1.improved} and Corollary~\ref{cor.0<q<1.improved}
are invariant under duality (of planar graphs),
which takes $v \mapsto q/v$ and interchanges loops and bridges.
In particular, the interval $(-1 - \sqrt{1-q}, -1 + \sqrt{1-q})$
is mapped onto itself under duality, with the endpoints interchanged.

\proofof{Theorem~\ref{thm.0<q<1.improved}}
The proof is by induction on the number of edges in $G$.
If $G$ has no edges, then $c=n$ and $Z_G(\qvbf) = q^n > 0$,
so the claim is obviously true.
Now suppose that $G$ has $m$ edges,
and assume that the result holds for all graphs having fewer than $m$ edges.
We now consider five cases:

(a)  If $G$ has a loop $e$, then
$Z_G(\qvbf) = (1+v_e) Z_{G \setminus e}(q, {\bf v}_{\neq e})$
by \reff{eq.delcon.loop}.
By hypothesis we have $1 + v_e > 0$;
and $G \setminus e$ has the same numbers of components and vertices
as $G$ does.
This proves that $Z_G(\qvbf)$ is nonzero with the desired sign.

(b)  If $G$ has a bridge $e$, then
$Z_G(\qvbf) = (q+v_e) Z_{G/e}(q, {\bf v}_{\neq e})$
by \reff{eq.delcon.bridge}.
By hypothesis we have $q + v_e < 0$;
and $G/e$ has the same numbers of components as $G$
but one less vertex.
This proves once again that $Z_G(\qvbf)$ is nonzero with the desired sign.

We can henceforth assume that $G$ has no loops or bridges.

(c) If $G$ has a pair $e_1,e_2$ of parallel edges, then we apply
the parallel-reduction formula \reff{eq.parallel2} to it.
By hypothesis, both $1+v_1$ and $1+v_2$ lie in the interval
$(-\sqrt{1-q},\sqrt{1-q})$.
Therefore $(1+v_1)(1+v_2)$ lies in the interval $(-(1-q),1-q)$,
and hence $v_{\rm eff}  \equiv (1+v_1)(1+v_2) - 1$
satisfies
\be
   v_{\rm eff}
   \;\in\; \Bigl( -1-(1-q) ,\, -1+(1-q) \Bigr)
   \;\subset\; \Bigl( -1-\sqrt{1-q} ,\, -q \Bigr)
   \;.
 \label{eq.v1parv2}
\ee
In the graph $G \setminus e_2$,
the edge $e_1$ is either a normal edge or a bridge; it cannot be a loop.
The new weight \reff{eq.v1parv2} satisfies the hypotheses
for both normal edges and bridges, so we may apply the inductive hypothesis
to $G \setminus e_2$.
Since $G \setminus e_2$ has the same number of vertices and components
as $G$ does, we are done.

(d) If $G$ has a pair $e_1,e_2$ of series edges {\em in the wide sense}\/,
then we apply the series-reduction formula \reff{eq.series2} to it.
By hypothesis, both $1+q/v_1$ and $1+q/v_2$ lie in the interval
$(-\sqrt{1-q},\sqrt{1-q})$.
Therefore $(1+q/v_1)(1+q/v_2)$ lies in the interval $(-(1-q),1-q)$,
hence
\be
   {q \over v_{\rm eff}}
   \;\equiv\; \biggl( 1 + {q \over v_1} \biggr)
              \biggl( 1 + {q \over v_2} \biggr) \,-\, 1
   \;\in\; (-2+q,-q)
   \;,
\ee
hence
\be
   v_{\rm eff}  \;\in\;  \Bigl( -1 ,\, -{q \over 2-q} \Bigr)
   \;\subset\; \Bigl( -1 ,\, -1+\sqrt{1-q} \Bigr)
   \;.
 \label{eq.v1serv2}
\ee
In the graph $G/e_2$,
the edge $e_1$ is either a normal edge or a loop; it cannot be a bridge.
The new weight \reff{eq.v1serv2} satisfies the hypotheses
for both normal edges and loops, so we may apply the inductive hypothesis
to $G/e_2$.
Now, $G/e_2$ has the same number of components as $G$ but one less vertex.
On the other hand, since $v_1, v_2 < -1+\sqrt{1-q} < -q/2$,
the prefactor $q+v_1+v_2$ in \reff{eq.series2} is negative.
This gives the correct sign.

(e) If $G$ has neither parallel edges nor
series edges {\em in the wide sense}\/,
then pick any edge $e$ and apply the deletion-contraction identity
\reff{eq.delcon} to it.
We see that $G \setminus e$ has the same number of vertices and components
as $G$ does (because $e$ is not a bridge);
and all edges of $G \setminus e$ are normal
(because $e$ does not belong to a wide-sense series pair in $G$,
 so a bridge cannot be formed by deletion).
Therefore, we can apply the inductive hypothesis to $G \setminus e$,
and the contribution has the correct sign.
Likewise, we see that $G/e$ has the same number of components as $G$
but one less vertex (because $e$ is not a loop);
and all edges of $G/e$ are normal
(because $e$ does not belong to a parallel pair in $G$,
 so a loop cannot be formed by contraction).
Therefore, we can apply the inductive hypothesis to $G/e$;
and since $v_e < 0$, the contribution again has the correct sign.
\qed

\bigskip
\bigskip

The following examples show that Theorem \ref{thm.0<q<1.improved} and 
Corollary~\ref{cor.0<q<1.improved} 
are in some sense 
best possible.
If $G$ is any tree,
we have $Z_G(q,v) = q (q+v)^{n-1}$,
so that there are roots at $v=-q$.
If $G$ has one vertex and $k$ loops, then $Z_G(q,v)=(1+v)^k$, 
so that there are roots at $v=-1$.
If $G$ is a cycle of length two, then  
we have $Z_G(q,v) = q (q+2v+v^2)$,
so that there are roots at $v = -1 \pm \sqrt{1-q}$.
We will see in Section \ref{subsec.0<q<1.block}, however, that 
Theorem \ref{thm.0<q<1.improved} can be improved if we add a hypothesis 
on the minimum number of edges in a block of $G$.

We next show that,
in the situation of Corollary~\ref{cor.0<q<1.improved}(a),
we can go farther and control derivatives with respect to $q$:

\begin{theorem}
   \label{thm.0<q<1.derivs}
Let $G$ be a loopless graph with $n$ vertices and $c$ components,
and let $q \in (0,1)$.
Suppose that $-1 - \sqrt{1-q} < v_e < -q$ for all $e \in E$.
Then for $0 \le \ell \le n-c$ we have
\be
   (-1)^{n-c-\ell} \, {d^\ell \over dq^\ell} \,
       \biggl( {Z_G(q,\bv) \over q^c} \biggr)
     \;>\; 0
   \;.
\ee
\end{theorem}

\proof
The proof is by induction on $|E|$.
If $G$ has no edges, then $c=n$ and hence $\ell=0$;
and since $Z_G(q,\bv)/q^c \equiv 1$, the result holds.

Now suppose that $G$ has $m$ edges,
and assume that the result holds for all graphs having fewer than $m$ edges.
We shall consider three cases:

(i)  If $G$ is a forest, then
then $Z_G(\qvbf)/q^c = \prod_{e \in E} (q+v_e)$.
Since $q+v_e < 0$ for all $e$,
this product and all its derivatives have the claimed sign.

We can henceforth assume that $G$ is not a forest.

(ii)  If $G$ has somewhere a pair $e_1,e_2$ of parallel edges,
pick some such pair and apply to it the parallel-reduction formula
\reff{eq.parallel2}, differentiated $\ell$ times with respect to $q$.
Since $(v_1,v_2) \mapsto (1+v_1)(1+v_2) - 1$
maps the interval $(-1 - \sqrt{1-q}, -q)$ into itself,
we can apply the inductive hypothesis to $G \setminus e_2$;
and the result has the right sign, since $G \setminus e_2$
has the same number of vertices and connected components as $G$ does.

(iii) If $G$ has no pair of parallel edges,
pick any non-bridge edge $e$ and apply to it
the deletion-contraction identity \reff{eq.delcon},
differentiated $\ell$ times with respect to $q$.
By the inductive hypothesis,
the first term has the right sign (strictly) because $G \setminus e$
has the same number of vertices and connected components as $G$ does
(since $e$ is not a bridge).
As for $G/e$, it has the same number of connected components as $G$
but one less vertex (since $e$ is not a loop).
Moreover, $G/e$ is loopless (since $G$ has no parallel edges).
If $\ell=n-c$,
then $(\partial^\ell/\partial q^\ell) (Z_{G/e}/q^c) = 0$
because $Z_{G/e}/q^c$ is a polynomial in $q$ of degree $n-1-c$
(note that $c \le n-1$ because $G$ has a non-loop edge $e$).
If $\ell < n-c$, we will be able to apply the inductive hypothesis to $G/e$;
and since $v_e < 0$, this term has the right sign as well
(strictly, though we do not need this).
\qed

Please note the strategy behind the proof of Theorem~\ref{thm.0<q<1.derivs}:
since the deletion-contraction and parallel-reduction formulae
do not involve $q$, they commute with differentiation with respect to $q$.
The series-reduction formula \reff{eq.series2}, by contrast,
involves $q$ both in the prefactor and (what seems to be worse)
in the argument $v_{\rm eff} = v_1 v_2 / (q+v_1+v_2)$;
we do not see how to handle the derivatives with respect to $q$.
It is for this reason that we limited ourselves to a situation
in which we could avoid the use of series reduction.
We do not know whether this restriction is really necessary.

\begin{corollary}
   \label{thm.0<q<1.derivs.chrompoly}
Let $G$ be a loopless graph with $n$ vertices and $c$ components,
and let $P_G(q)$ be its chromatic polynomial. Then for all $\ell \ge 0$,
\be
     (-1)^{n-c-\ell} \, {d^\ell \over dq^\ell} \, {P_G(q) \over q^c}
     \Biggr| _{q=1}
     \;\ge\;  0   \;.
 \label{eq.thm.0<q<1.derivs.chrompoly}
\ee
\end{corollary}

\proof
For $0 \le \ell \le n-c$,
put $v_e = -1$ for all $e$ in Theorem~\ref{thm.0<q<1.derivs}
and let $q \uparrow 1$.
For $\ell > n-c$, \reff{eq.thm.0<q<1.derivs.chrompoly} of course vanishes.
\qed

\medskip\par\noindent
{\bf Remark.}  Since the numbers \reff{eq.thm.0<q<1.derivs.chrompoly}
are nonnegative integers, it would be nice to find a combinatorial
interpretation for them.
Since $P_G(q)/q^c$ factorizes over blocks,
it suffices to do this for 2-connected graphs $G$.
For $\ell = 1$ the following characterization is known:
For any connected graph $G$ on $n$ vertices
and any edge $e = ij$ of $G$,
the quantity $[\partial T_G(x,y)/\partial x](0,0) = (-1)^n P'_G(1)$
counts the acyclic orientations of $G$
in which $i$ is the unique source and $j$ is the unique sink
\cite[Theorem 7.2]{Greene_83}  \cite[Exercise 6.35]{Brylawski_92}
\cite{Gebhard_00}.
If $G$ is bridgeless, then $(-1)^n P'_G(1)$
also equals half the number of totally cyclic orientations of $G$
(i.e.\ orientations in which every edge of $G$ belongs to some
 directed cycle)
in which every directed cycle uses $e$
\cite[Theorem 8.2]{Greene_83}
\cite[Proposition 6.2.12 and Example 6.3.29]{Brylawski_92}.
\qed

\bigskip

The matroidal version of Theorem~\ref{thm.0<q<1.improved}
is proven by an identical argument:

\begin{theorem}
   \label{thm.0<q<1.improved.matroid}
Let $M$ be a matroid with ground set $E$ and rank $r(M)$,
and let $q \in (0,1)$.
Suppose that:
\begin{itemize}
   \item[(i)]  $v_e > -1$ for every loop $e$;
   \item[(ii)]  $v_e < -q$ for every coloop $e$; and
   \item[(iii)]  $-1 - \sqrt{1-q} < v_e < -1 + \sqrt{1-q}$
       for every normal (i.e., non-loop non-coloop) element $e$.
\end{itemize}
Then $(-1)^{r(M)} \Ztilde_M(\qvbf) > 0$.
\end{theorem}

The matroidal analogue of Corollary~\ref{cor.0<q<1.improved} 
is obvious, and we refrain from stating it explicitly.

Finally, we have the following matroidal version of
Theorem~\ref{thm.0<q<1.derivs}:

\begin{theorem}
   \label{thm.0<q<1.derivs.matroid}
Let $M$ be a loopless matroid with ground set $E$ and rank $r(M)$,
and let $q \in (0,1)$.
Suppose that $-1 - \sqrt{1-q} < v_e < -q$ for all $e \in E$.
Then for $0 \le \ell \le r(M)$ we have
\be
   (-1)^{r(M)-\ell} \, {d^\ell \over dq^\ell} \,
       \Bigl( q^{r(M)} \Ztilde_M(q,\bv) \Bigr)
     \;>\; 0
   \;.
\ee
\end{theorem}

%
%
%
%
%

\section{An abstract theorem}  \label{sec.abstract}


In this section we prove an abstract result
that we shall subsequently use in two ways:
in Section~\ref{subsec.0<q<1.block} we will use it to strengthen
Corollary~\ref{cor.0<q<1.improved} by considering the block structure of $G$,
and in Section~\ref{sec.q>1} we will use it
to obtain zero-free regions when $q \in (1,32/27]$.

Since our ``graphs'' allow loops and multiple edges,
let us be completely precise about what we mean by ``blocks''.
We say that a graph $G=(V,E)$ is {\em separable}\/
if there exist graphs $G_1=(V_1,E_1)$ and $G_2=(V_2,E_2)$
such that $G = G_1 \cup G_2$, $G \neq G_1$, $G \neq G_2$,
$E_1 \cap E_2 = \emptyset$ and $|V_1 \cap V_2| \le 1$.
A {\em block}\/ of $G$ is a maximal non-separable subgraph of $G$.
We say that a graph $G=(V,E)$ is {\em $k$-connected}\/ ($k \ge 2$)
in case it has at least $k+1$ vertices
and $G \setminus U$ is connected for all $U\subseteq V$ with $|U|<k$.
(Thus, a graph is $k$-connected if and only if
 its underlying simple graph is $k$-connected.)
Let us remark that a graph is non-separable
if and only if it is either 2-connected and loopless
or else is $K_1$ (a single vertex with no edges),
$C_1$ (a single vertex with a loop)
or $K_2^{(m)}$ (a pair of vertices connected by $m$ parallel edges, $m \ge 1$).
Equivalently, a graph $G \neq K_1, K_2$ is non-separable
if and only if it has no isolated vertices
and every pair of distinct edges belongs to a cycle.
Finally (and most importantly),
a graph $G \neq K_1$ is non-separable
if and only if it has no isolated vertices
and its cycle matroid $M(G)$ is $2$-connected.

We recall that a graph $H$ is a {\em minor}\/ of a graph $G$
(written $H \preceq G$) in case $H$ can be obtained from $G$
by a sequence (possibly empty) of deletions of edges,
contractions of edges, and deletions of isolated vertices.
Note, in particular, that any subgraph of $G$ is a minor of $G$.
Note also that parallel and series reduction lead to minors,
because they are special cases of edge deletion and contraction, respectively.
A class $\scrg$ of graphs is a {\em minor-closed class}\/
in case $G \in \scrg$ and $H \preceq G$ imply $H \in \scrg$.

\begin{theorem}
   \label{thm.0<q<1.block.abstract}
Let $m \ge 2$, $\gamma\in \{0,1\}$ and $q > 0$.
Let $\scrg$ be a minor-closed class of graphs.
Suppose that $\scrv \subset \R$ satisfies the following hypotheses:
\begin{itemize}
   \item[(a)]  $\scrv \subseteq (-2,-q/2)$
   \item[(b)]  $\scrv \parallel \scrv \subseteq \scrv$ 
[Here $\parallel$ denotes parallel connection, as defined after
       \reff{eq.parallel2}.]
   \item[(c)]  $\scrv \series \scrv \subseteq \scrv$  [Here $\series$ denotes series connection, as defined after
       \reff{eq.series2a}.]
   \item[(d)]  $(-1)^{|V|-1+\gamma} Z_G(\qvbf) > 0$ whenever $G=(V,E) \in \scrg$ is a
       non-separable graph with {\em exactly}\/ $m$ edges,
       and $v_e \in \scrv$ for all $e \in E$.
\end{itemize}
Then $(-1)^{n+c+\gamma b} Z_G(\qvbf) > 0$
whenever $G=(V,E) \in \scrg$ is a graph with $n$ vertices, $c$ components
and $b$ blocks,
in which each block contains {\em at least}\/ $m$ edges,
and $v_e \in \scrv$ for all $e \in E$.
\end{theorem}

Before proving Theorem~\ref{thm.0<q<1.block.abstract},
let us make a few simple observations about the hypotheses and the conclusion:

1)  If the set $\scrv \subset \R$ satisfies hypotheses (a)--(d)
for a given $m$, then it also satisfies those hypotheses
for all larger $m$;
this is not {\em a priori}\/ obvious for hypothesis (d),
but it is part of the conclusion of the theorem.

2) The conditions (a)--(c) on $\scrv$ are invariant under the duality map
$v \mapsto q/v$, and the class of {\em connected planar}\/ graphs
in which each block contains exactly (resp.\ at least) $m$ edges
is also invariant under duality.

3)  In the presence of hypothesis (b),
hypothesis (a) is equivalent to the weaker condition
$\scrv \subseteq (-\infty,-q/2)$,
since $v \parallel v \ge 0 > -q/2$ whenever $v \le -2$.
Indeed, the condition $\scrv \subseteq (-\infty,-q/2)$
is all that is actually used in the proof of
Theorem~\ref{thm.0<q<1.block.abstract}.
We have stated hypothesis (a) in the stronger form
in order to make manifest the duality-invariance.

4)  The proof of Theorem~\ref{thm.0<q<1.block.abstract} uses only $q > 0$,
but we will show in Corollary~\ref{propdiamond}
that hypotheses (a)--(c) can be satisfied (with $\scrv \neq \emptyset$)
only if $q \le 32/27$ and $\scrv$ is contained in a particular interval
$I_\Diamond(q)$.
In addition, we will show in Proposition~\ref{gamma}
that hypotheses (b) and (d) can be satisfied
(with $\scrv \neq \emptyset$ and $\scrg \supseteq$ series-parallel graphs)
only if $\gamma=0$ and  $q< 1$, or $\gamma=1$ and $q > 1$.  
Finally, in Corollary~\ref{cor.diamond2a} we will show that
the {\em conclusion}\/ of Theorem~\ref{thm.0<q<1.block.abstract} can hold
(with $\scrv \neq \emptyset$ and $\scrg \supseteq$ series-parallel graphs)
only if either $\gamma=0$, $q < 1$ and $\scrv \subseteq I_\Diamond(q)$
or else $\gamma=1$, $1 < q \le 32/27$ and $\scrv \subseteq I_\Diamond(q)$.

5) We shall be principally interested in the case when $\scrg = $ all graphs,
but we have stated Theorem~\ref{thm.0<q<1.block.abstract}
for an arbitrary minor-closed class because it is no more difficult
to prove, and other minor-closed classes
(e.g.\ planar graphs, series-parallel graphs)
may be of interest.

\bigskip

The proof of Theorem~\ref{thm.0<q<1.block.abstract}
will of course be based on deletion-contraction
(together with parallel and series reduction);
but it will be slightly more delicate than the proofs
of Theorems~\ref{thm.0<q<1.improved} and \ref{thm.0<q<1.derivs},
because in order to apply the inductive hypothesis,
we will need to find an element $e$ for which
{\em both}\/ $G \setminus e$ and $G/e$ are non-separable.
(In Theorems~\ref{thm.0<q<1.improved} and \ref{thm.0<q<1.derivs}
 we needed only to maintain connectedness, not non-separability.)
A sufficient condition for this is provided by
the following graph-theoretic result:

\begin{proposition}
   \label{prop.combinatorial.graph}
Let $G$ be a simple 2-connected graph
with at most one vertex of degree 2.
Then there exists an edge $e$
such that both $G \setminus e$ and $G/e$ are 2-connected.
\end{proposition}

\proof
If $G$ is 3-connected, then it is easy to see that
$G \setminus e$ and  $G/e$ are both 2-connected for all $e\in E$.
Thus we may suppose that there exists
$U = \{u_1,u_2\} \subseteq V$ such that $G \setminus U$ is disconnected.
Fix a vertex $x_0 \in V$ such that all vertices of $G$ other than $x_0$
have degree $\ge 3$.
Choose a pair $(U,H)$ such that $U =  \{u_1,u_2\} \subseteq V$,
$G \setminus U$ is disconnected, $H$ is a component of $G \setminus U$,
$x_0\not\in V(H)$, and $|V(H)|$ is as small as possible
consistent with these constraints.
The fact that $G$ is simple and $d_G(x)\geq 3$ for all $x\in V(H)$
implies that $E(H)\neq \emptyset$.

Let $H_1$ (resp.\ $H_2$) be the subgraph of $G$ induced by $V(H)\cup U$
[resp.\ by $V \setminus V(H)$],
and let $H'_1$ (resp.\ $H'_2$) be the graph obtained
from $H_1$ (resp.\ $H_2$) by adding the edge $u_1 u_2$
if it is not already in $G$.
The minimality of $|V(H)|$ implies that $H'_1$ is 3-connected:
for if $H'_1$ had a 2-vertex cut $U'$,
then we must have $U' \neq U$
(since $H'_1 \setminus U = H$ is connected)
and it is not hard to see that some component of $H'_1 \setminus U'$
(indeed, any component disjoint from $U \setminus U'$)
is also a component of $G \setminus U'$ that is strictly contained in $H$.

Therefore, $H'_1 \setminus e$ and $H'_1/e$ are both 2-connected
for all $e\in E(H'_1)$.
On the other hand, $H'_2$ is also 2-connected.
So choose any $e \in E(H'_1)$, $e \neq u_1 u_2$,
glue $H'_1 \setminus e$ (resp.\ $H'_1/e$) onto $H'_2$ along $\{u_1,u_2\}$,
and delete the edge $u_1 u_2$;
this operation (2-sum) preserves 2-connectivity
and yields $G \setminus e$ (resp.\ $G/e$).
\qed

\medskip\par\noindent
{\bf Remarks.}
1.  We will actually need here only the weaker version of
Proposition~\ref{prop.combinatorial.graph} in which
``at most one vertex of degree 2'' is replaced by
``no vertices of degree 2'', i.e.\ $G$ has minimum degree $\ge 3$.

2. Please note that, by definition, $G$ is simple $\iff$ $G$ is loopless
and has no parallel edges;
and that, for loopless graphs,
$G$ has no vertex of degree 2 $\iff$ $G$ has no pair of series edges
in the narrow sense.
These trivial facts will be used in step (iii) of the proof of
Theorem~\ref{thm.0<q<1.block.abstract}.

3.  It is natural to ask, for arbitrary $k \ge 2$,
how large a minimum degree is needed in a $k$-connected simple graph
in order that there exist an edge $e$ such that both
$G \setminus e$ and $G/e$ are $k$-connected.
For $k=2$, Proposition~\ref{prop.combinatorial.graph}
gives the optimal answer: minimum degree at least 3.
For $k \ge 3$, a {\em sufficient}\/ condition is
minimum degree $\ge \lceil (3k-1)/2 \rceil$:
this follows from the result of Chartrand, Kaugars and Lick \cite{CKL}
that every $k$-connected simple graph
of minimum degree at least $\lceil (3k-1)/2 \rceil$
has a vertex $x$ such that $G \setminus x$ is $k$-connected;
then any edge $e$ incident on $x$ will do.
This gives the optimal answer also for $k=3$:
minimum degree at least 4.
For $k \ge 4$ the optimal result is apparently not known.
Note, however, that for any $k\geq 4$ there exist $k$-connected graphs of
minimum degree $\lfloor 5k/4 \rfloor - 1$
(but no higher) with no edges $e$ such that $G/e$ is $k$-connected
\cite[p.~16]{Egawa_91} \cite[p.~97]{Kriesell_01}.

\proofof{Theorem~\ref{thm.0<q<1.block.abstract}}
Since $Z_G(\qvbf)$ factorizes over blocks
(modulo a factor $q>0$ in the case of a cut vertex)
and the quantity $n-c+\gamma b$ is additive over blocks,
and every block of $G$ is a minor of $G$ (hence belongs to $\scrg$),
it suffices to prove Theorem~\ref{thm.0<q<1.block.abstract}
for non-separable graphs $G \in \scrg$.

The proof is by induction on $|E|$.
The base case is $|E| = m$, which holds by hypothesis (d).

Assume now that $|E| > m$.  We consider three cases:

(i) If $G$ has somewhere a pair $e_1,e_2$ of parallel edges,
pick some such pair and apply the parallel-reduction formula
\reff{eq.parallel2} to it.
Since $\scrv \parallel \scrv \subseteq \scrv$,
and $G \setminus e_2$ is non-separable and has at least $m$ edges
(and belongs to $\scrg$),
we can apply the inductive hypothesis to $G \setminus e_2$;
and the result has the right sign, since $G \setminus e_2$
has the same number of vertices as $G$ does.

(ii)  If $G$ has somewhere a pair $e_1,e_2$ of series edges
(in either the narrow sense or the wide sense, it doesn't matter),
pick some such pair and apply the series-reduction formula
\reff{eq.series2} to it.
Since $\scrv \subseteq (-\infty,-q/2)$,
the prefactor $q + v_{e_1} + v_{e_2}$ is $<0$;
furthermore, since $\scrv \series \scrv \subseteq \scrv$,
and $G / e_2$ is non-separable and has at least $m$ edges
(and belongs to $\scrg$),
we can apply the inductive hypothesis to $G / e_2$;
and the result has the right sign, since $G / e_2$
has one less vertex than $G$ does.

(iii) If $G$ has neither a pair of parallel edges
nor a pair of series edges, then $G$ is simple and has no degree-2 vertices,
so by Proposition~\ref{prop.combinatorial.graph}
there exists $e \in E$ such that both $G \setminus e$ and $G/e$
are 2-connected (and hence non-separable).
So we can use the deletion-contraction identity on $e$
and apply the inductive hypothesis to both $G \setminus e$ and $G/e$
(which belong to $\scrg$).
The result has the right sign,
because $G \setminus e$ (resp.\ $G/e$) has $|V|$ (resp.\ $|V|-1$) vertices
and $v_e < 0$.
\qed

\begin{proposition}
   \label{gamma}
Fix $m \ge 2$, $\gamma\in\{0,1\}$ and $q > 0$.
Suppose that $\scrv \subset \R$ ($\scrv \neq \emptyset$)
satisfies hypotheses (b) and (d) of Theorem~\ref{thm.0<q<1.block.abstract}
for a class $\scrg \supseteq$ series-parallel graphs.
Then either
\begin{itemize}
   \item[(a)]  $\gamma=0$, $q<1$ and $\scrv \subseteq (-2,0)$, or
   \item[(b)]  $\gamma=1$, $q > 1$ and $\scrv \subseteq (-2,0)$.
\end{itemize}
\end{proposition}

\proof
We begin with the trivial observation that $Z_G(q,\bv) > 0$
whenever $q > 0$ and $v_e \ge 0$ for all $e$.  It follows that:

(i) If $\scrv \cap [0,\infty) \neq \emptyset$,
we can obtain a counterexample to hypothesis (d)
--- no matter what the values of $m$, $\gamma$ and $q$ ---
by taking $G=(V,E)$ to be any 2-connected graph in $\scrg$
with exactly $m$ edges that has $|V|+\gamma$ even,
and then taking all $v_e \in \scrv \cap [0,\infty)$.

(ii) If $\scrv \cap (-\infty,-2]  \neq \emptyset$,
then closure under parallel connection [hypothesis (b)]
implies that $\scrv \cap [0,\infty) \neq \emptyset$,
so we are reduced to case (i).

We have therefore proven that $\scrv \subseteq (-2,0)$.

(iii) Since $\scrv \cap (-2,0) \neq \emptyset$,
it follows by repeated parallel connection (of even order) that $\scrv$
contains, for every $\epsilon > 0$, a point $v_\star \in [-1,-1+\epsilon]$.
Now consider the graph $G = K_2^{(m)}$,
for which we have
\be
   Z_{K_2^{(m)}}(q,v_1,\ldots,v_m)
   \;=\;
   q \left[ q-1 + \prod_{i=1}^m (1+v_i) \right]
   \;.
 \label{eq.prop.q>1.converse}
\ee
Taking $v_1 = \ldots = v_m = v_\star$
and letting $\epsilon \downarrow 0$,
we see that $(-1)^{1+\gamma} Z_{K_2^{(m)}} > 0$ requires
either $\gamma=0$ and $q<1$, or $\gamma=1$ and $q \geq 1$.

To rule out $q=1$, we observe that
$Z_G(1,\bv) = \prod\limits_{e \in E} (1+v_e)$.  Then:

(iv) If $\scrv \cap [-1,0) \neq \emptyset$,
we can obtain a counterexample to hypothesis (d)
by taking $G=(V,E)$ to be any 2-connected graph in $\scrg$
with exactly $m$ edges and $|V|+\gamma$ even.

(v) If $\scrv \cap (-2,-1] \neq \emptyset$,
we can obtain a counterexample to hypothesis (d)
by taking $G=(V,E)$ to be any 2-connected graph in $\scrg$
with exactly $m$ edges and $|E|-|V|+\gamma$ even.
\qed

{\bf Remark.}
It is obvious from the proof that this result
holds for classes $\scrg$ much smaller than all series-parallel graphs.

\bigskip

We conclude this section by giving the matroidal analogue of
Theorem~\ref{thm.0<q<1.block.abstract}.
We shall be brief, because the proofs are nearly identical
to the proofs for graphs;  we shall merely point out the differences.

Recall first that a matroid $N$ is a {\em minor}\/ of a matroid $M$
(written $N \preceq M$) in case $N$ can be obtained from $M$
by a sequence (possibly empty) of deletions or contractions of elements.
In particular, parallel and series reduction lead to minors,
because they are special cases of deletion and contraction, respectively.
A class $\scrm$ of matroids is a {\em minor-closed class}\/
in case $M \in \scrm$ and $N \preceq M$ imply $N \in \scrm$.
(In particular, {\em graphic}\/ matroids form a minor-closed class.)

\begin{theorem}
   \label{thm.0<q<1.block.abstract.matroid}
Let $m \ge 2$, $\gamma\in \{0,1\}$ and $q > 0$.
Let $\scrm$ be a minor-closed class of matroids.
Suppose that $\scrv \subset \R$ satisfies the following hypotheses:
\begin{itemize}
   \item[(a)]  $\scrv \subseteq (-2,-q/2)$
   \item[(b)]  $\scrv \parallel \scrv \subseteq \scrv$ (parallel connection)
   \item[(c)]  $\scrv \series \scrv \subseteq \scrv$ (series connection)
   \item[(d)]  $(-1)^{r(M)+\gamma} \Ztilde_M(\qvbf) > 0$ whenever
       $M \in \scrm$ is a 2-connected matroid with ground set $E$ and
       {\em exactly}\/ $m$ elements (i.e., $|E|=m$)
       and $v_e \in \scrv$ for all $e \in E$.
\end{itemize}
Then $(-1)^{r(M)+\gamma b} \Ztilde_M(\qvbf) > 0$
whenever $M \in \scrm$ is a matroid with ground set $E$ with $b$ 2-connected
components,
in which each 2-connected component contains {\em at least}\/ $m$ elements,
and $v_e \in \scrv$ for all $e \in E$.
\end{theorem}

\noindent
Theorem~\ref{thm.0<q<1.block.abstract} is simply the special case
of Theorem~\ref{thm.0<q<1.block.abstract.matroid}
in which $\scrm$ is a minor-closed class of {\em graphic}\/ matroids.

The key to the proof of Theorem~\ref{thm.0<q<1.block.abstract.matroid}
is the following matroidal analogue to
Proposition~\ref{prop.combinatorial.graph},
which was proven by Oxley \cite[Proposition~3.5]{Oxley_84}:

\begin{proposition}[Oxley \protect\cite{Oxley_84}]
   \label{prop.Oxley}
Let $M$ be a 2-connected matroid having at least 2 elements,
and let $d_M(k)$ [resp.\ $d^*_M(k)$] be the number of $k$-element
circuits (resp.\ cocircuits) in $M$.
If $d_M(2) + d^*_M(2) \le 1$,
then there exists an element $e\in E$
for which both $M \setminus e$ and $M/e$ are 2-connected.
\end{proposition}

Once again, we shall need only the special case of this result
for $d_M(2) = d^*_M(2) = 0$,
i.e.\ when there are no 2-element circuits (= pairs of parallel elements)
or 2-element cocircuits (= pairs of series elements).
The proof of Theorem~\ref{thm.0<q<1.block.abstract.matroid}
is then identical to that of Theorem~\ref{thm.0<q<1.block.abstract},
but using Proposition~\ref{prop.Oxley} in place of
Proposition~\ref{prop.combinatorial.graph}.
Here we are obliged to understand ``series elements'' in the wide sense,
since this is the only sense that makes sense for matroids.

\section{The interval {\boldmath $q\in (0,1)$} revisited}
   \label{subsec.0<q<1.block} \label{sec4.2}

We believe that Corollary~\ref{cor.0<q<1.improved}(c) is
the first of an infinite family of results
giving successively larger zero-free regions
under successively stronger hypotheses on the size of the blocks
that can appear in $G$.
Stating that $G$ is loopless and bridgeless is equivalent
to saying that each block of $G$
(other than possible isolated vertices) contains at least two edges.
Furthermore, the extremal graph for Corollary~\ref{cor.0<q<1.improved}(c)
is the unique block with exactly two edges, namely $C_2 = K_2^{(2)}$.
Similarly, there will be theorems stating that if each block of $G$
contains at least $m$ edges, then $(-1)^{n+c} Z_G(\qvbf) > 0$
whenever all the $v_e$ lie in a particular 
(maximal) interval $(v_m^-(q), v_m^+(q))$.
We {\em conjecture}\/ that $v_m^+(q)$ [resp.\ $v_m^-(q)$]
can be chosen to be 
{\em strictly}\/ increasing (resp.\ decreasing) in $m$.

We can use 
Theorem~\ref{thm.0<q<1.block.abstract}
to determine an interval
$\scrv = (v_m^-(q), v_m^+(q))$ for the cases $m=2,3,4$
and $\scrg = $ all graphs.
Let us begin with a simple result that is optimal for $m=2,4$
but not for $m=3$.

\begin{corollary}
   \label{cor.0<q<1.block34}
Let $2 \le m \le 4$,
and let $G$ be a graph with $n$ vertices and $c$ components,
in which each block contains at least $m$ edges.
Let $0 < q < 1$, and suppose that
\be
   -[1 + (1-q)^{1/m}]  \;<\;  v_e  \;<\;  - \, {q \over 1 + (1-q)^{1/m}}
 \label{eq.cor.0<q<1.block34}
\ee
for all $e \in E$.
Then $(-1)^{n+c} Z_G(\qvbf) > 0$.
\end{corollary}

\noindent
The case $m=2$ is of course just Corollary~\ref{cor.0<q<1.improved}(c);
the new cases here are $m=3$ and $m=4$.
Let us remark that the intervals in \reff{eq.cor.0<q<1.block34}
are invariant under the duality map $v \mapsto q/v$.

Let us begin by working out the special cases in which $G$ is either
an $m$-cocycle $K_2^{(m)}$ (two vertices connected by $m$ parallel edges)
or an $m$-cycle $C_m$.
This calculation works for all $m$:

\begin{lemma}
   \label{lemma.0<q<1.K2mCm}
Let $m \ge 1$ and $0 < q < 1$.  Then:
\begin{itemize}
   \item[(a)] If $G=K_2^{(m)}$ and
      $- [1 + (1-q)^{1/m}] < v_e < -[1 - (1-q)^{1/m}]$ for all $e$,
      then $Z_{K_2^{(m)}}(\qvbf) < 0$.
Conversely, $Z_{K_2^{(m)}}(q,v) = 0$ when
\begin{subeqnarray}
   v & = & -[1 - (1-q)^{1/m}] \quad\hbox{(all $m$)}  \\
   v & = & -[1 + (1-q)^{1/m}] \quad\hbox{(even $m$)}
\end{subeqnarray}
   \item[(b)] If $G=C_m$ and
      $-q/[1 - (1-q)^{1/m}] < v_e < -q/[1 + (1-q)^{1/m}]$ for all $e$,
      then $(-1)^m Z_{C_m}(\qvbf) < 0$.
Conversely, $Z_{C_m}(q,v) = 0$ when
\begin{subeqnarray}
   v & = & -q/[1 - (1-q)^{1/m}] \quad\hbox{(all $m$)}  \\
   v & = & -q/[1 + (1-q)^{1/m}] \quad\hbox{(even $m$)}
\end{subeqnarray}
\end{itemize}
\end{lemma}

\proof
We have
\be
   Z_{K_2^{(m)}}(\qvbf)  \;=\;
   q \left[ q-1 + \prod\limits_{i=1}^m (1+v_i) \right]
   \;.
 \label{eq.Z.K2m}
\ee
If $0 < q < 1$ and $|1+v_i| < (1-q)^{1/m}$ for all $i$,
we obviously have $Z_{K_2^{(m)}}(\qvbf) < 0$.
The converse claims follow easily from \reff{eq.Z.K2m}.
This proves (a).

Part (b) then follows by using the duality relation \reff{eq.duality.matroid},
noting that if $G=(V,E)$ is a planar graph and $G^*=(V^*,E^*)$ is its dual,
then $(-1)^{|V|-1} = (-1)^{r(M(G))} = (-1)^{E - r(M^*(G))}
      = (-1)^{|E^*|} (-1)^{|V^*|-1}$
and the prefactor $\prod_{e \in E} v_e$ in \reff{eq.duality.matroid}
has sign $(-1)^{E} = (-1)^{|E^*|}$.
\qed

{\bf Remarks.}
1.  The foregoing argument also shows that $Z_{K_2^{(m)}}(\qvbf) \neq 0$
for {\em complex}\/ $\bv$ lying in the disc $|1+v_i| < (1-q)^{1/m}$
for all $i$.
Likewise, $Z_{C_m}(\qvbf) \neq 0$ for complex $\bv$
lying in the dual disc $|1+q/v_i| < (1-q)^{1/m}$,
i.e.\ the disc in complex $\bv$-space whose diameter is the interval
$\big(-q/[1 - (1-q)^{1/m}], \, -q/[1 + (1-q)^{1/m}] \big)$.

2.  It follows from the converse half of Lemma~\ref{lemma.0<q<1.K2mCm}
that, for $m=2$ and $m=4$,
the interval \reff{eq.cor.0<q<1.block34} is best possible
[and in general that, for even $m$, one cannot possibly do better
than \reff{eq.cor.0<q<1.block34}].
Indeed, this interval is best possible even in the univariate case.
We shall discuss the case $m=3$ after completing the proof of
Corollary~\ref{cor.0<q<1.block34}.
\qed

\bigskip

In the light of Lemma~\ref{lemma.0<q<1.K2mCm}, let us define the intervals
\begin{eqnarray}
   I_m^{\rm cocyc}
   & = &
   \Bigl( - [1 + (1-q)^{1/m}], \: - [1 - (1-q)^{1/m}]  \Bigr)
       \label{def.Imcocyc} \\[3mm]
   I_m^{\rm cyc}
   & = &
   \biggl( - \, {q \over 1 - (1-q)^{1/m}} ,\:
           - \, {q \over 1 + (1-q)^{1/m}}     \biggr)
       \label{def.Imcyc}
\end{eqnarray}
for arbitrary real $m \ge 1$.
Note that for $m \ge 2$ we have
\be
   - \, {q \over 1 - (1-q)^{1/m}}
   \;\le\;
   -[1 + (1-q)^{1/m}]
   \;<\;
   - \, {q \over 1 + (1-q)^{1/m}}
   \;\le\;
   -[1 - (1-q)^{1/m}]
   \;,
  \label{eq.Imineq}
\ee
while for $1 \le m \le 2$ we have the reverse inequality
\be
   -[1 + (1-q)^{1/m}]
   \;\le\;
   - \, {q \over 1 - (1-q)^{1/m}}
   \;<\;
   -[1 - (1-q)^{1/m}]
   \;\le\;
   - \, {q \over 1 + (1-q)^{1/m}}
   \;.
  \label{eq.Imineq.reverse}
\ee
In particular, for $m \ge 2$ the intersection
$I_m^{\rm cocyc} \cap I_m^{\rm cyc}$ is the self-dual interval
\be
   I_m  \;\equiv\;
   \biggl( - [1 + (1-q)^{1/m}], \: - \, {q \over 1 + (1-q)^{1/m}} \biggr)
   \;,
 \label{def.Im}
\ee
which is precisely the interval that arises (for $m=2,3,4$)
in Corollary~\ref{cor.0<q<1.block34}.
This interval has the following easily-verified properties:

\begin{lemma}
   \label{lemma.Im}
For any $q \in (0,1)$, the intervals
$I_m$ defined by \reff{def.Im} have the following properties:
\begin{itemize}
   \item[(a)]  $I_m$ is self-dual, i.e.\ it is invariant under $v \mapsto q/v$.
   \item[(b)]  If $m < m'$, then $I_m \subsetneq I_{m'}$.
   \item[(c)]  $I_4 \parallel I_4 \subseteq I_2$.
   \item[(d)]  $I_4 \series I_4 \subseteq I_2$.
   \item[(e)]  $\lim\limits_{m \to\infty} I_m = (-2, -q/2)$.
\end{itemize}
\end{lemma}

\proof
(a), (b) and (e) are obvious.
To prove (c), note that
$I_4 \subset I_4^{\rm cocyc}$ by \reff{eq.Imineq},
while $I_4^{\rm cocyc} \parallel I_4^{\rm cocyc} = I_2^{\rm cocyc} = I_2$.
Statement (d) follows from (c) and self-duality.
\qed

We are now ready to prove Corollary~\ref{cor.0<q<1.block34}:

\proofof{Corollary~\ref{cor.0<q<1.block34}}
We need only verify that the interval $\scrv = I_m$
defined in \reff{def.Im} satisfies hypotheses (a)--(d) of
Theorem~\ref{thm.0<q<1.block.abstract} with $\gamma=0$.
Hypothesis (a) follows from Lemma~\ref{lemma.Im}(b,e).
Since $2 \le m \le 4$,
hypothesis (b) follows from
$I_m \parallel I_m \subseteq I_4 \parallel I_4 \subseteq I_2 \subseteq I_m$
by Lemma~\ref{lemma.Im}(b,c).
Hypothesis (c) follows from hypothesis (b) by Lemma~\ref{lemma.Im}(a).

To prove hypothesis (d),
we must consider all non-separable graphs with $m$ edges.
For $m=2$ and $m=3$, the only such graphs are $m$-cocycles and $m$-cycles,
so the required statement follows from Lemma~\ref{lemma.0<q<1.K2mCm}
together with the observation \reff{eq.Imineq}.
For $m=4$ we must also consider the triangle with one double edge
(which can alternatively be thought of as the wheel $W_2$).
Applying parallel reduction to the double edge
and series reduction to the other pair of edges,
and using $I_4 \parallel I_4 \subseteq I_2$
and $I_4 \series I_4 \subseteq I_2$
from Lemma~\ref{lemma.Im}(c,d),
we reduce to the case of a 2-cycle (= 2-cocycle) with edge weights in $I_2$.
This proves hypothesis (d) for $m=4$.
\qed

So the interval \reff{eq.cor.0<q<1.block34} is optimal for $m=2,4$,
i.e.\ we have
\be
   \left\{  \begin{array}{l}
                v_m^-(q) \;=\; -[1 + (1-q)^{1/m}]  \\[2mm]
                v_m^+(q) \;=\; -q/[1 + (1-q)^{1/m}]
            \end{array}
   \right\}
   \quad\hbox{for }  m=2,4
 \label{def.v24}
\ee

For $m=3$, the interval \reff{eq.cor.0<q<1.block34} is {\em not}\/ optimal;
indeed, we suspect that there is no single optimal interval.
(That is, it may be possible, starting from an optimal interval,
to simultaneously increase or simultaneously decrease
both $v_3^-$ and $v_3^+$, yielding an incomparable optimal interval.)
To avoid these complications, let us restrict attention
to {\em self-dual}\/ intervals $\scrv$,
i.e.\ $\scrv = (v_m^-(q), v_m^+(q))$ with $v_m^-(q) = q/v_m^+(q)$.
For intervals of this kind we can state an optimal result for $m=3$:

\begin{corollary}
   \label{cor.0<q<1.block3optimal}
Let $G$ be a graph with $n$ vertices and $c$ components,
in which each block contains at least three edges.
Let $0 < q < 1$, and suppose that
\be
   {q \over v_3^+(q)}  \;<\;  v_e  \;<\; v_3^+(q)
\ee
for all $e \in E$, where $v_3^+(q)$ is the unique real root
of the cubic equation
\be
   v^3 \,+\, 3qv^2 \,+\, (q^2+2q)v \,+\, q^2  \;=\;  0 \;.
 \label{eq.cor.0<q<1.block3optimal.cubic}
\ee
Then $(-1)^{n+c} Z_G(\qvbf) > 0$.
\end{corollary}

\proof
First let us show that for $0 < q < 1$, the cubic equation
\reff{eq.cor.0<q<1.block3optimal.cubic} does indeed have a
single real root $v_3^+(q)$,
which lies between $-q$ and $-q/2$.
This is easy:  the derivative of the cubic
\reff{eq.cor.0<q<1.block3optimal.cubic},
namely $3v^2 + 6qv + (q^2+2q)$,
has discriminant $24(q^2-q) < 0$,
so the cubic \reff{eq.cor.0<q<1.block3optimal.cubic}
has strictly positive derivative on all of $\R$.
Moreover, the cubic \reff{eq.cor.0<q<1.block3optimal.cubic}
takes the value $q^3/8 > 0$ at $v=-q/2$
and the value $q^3 - q^2 < 0$ at $v=-q$,
so the unique real root must lie between $-q$ and $-q/2$.

Now consider $\scrv = (q/v_+,v_+)$ with $-2 < q/v_+ < -1 < -q < v_+ < -q/2$,
and let us try to satisfy
hypotheses (a)--(d) of Theorem~\ref{thm.0<q<1.block.abstract} with
$\gamma=0$.
We need to choose $v_+$ so that $\scrv \parallel \scrv \subseteq \scrv$
and $Z_{K_2^{(3)}}(q,v_1,v_2,v_3) < 0$ for $v_1,v_2,v_3 \in \scrv$.
If we succeed in doing this, then duality will guarantee that
$\scrv \series \scrv \subseteq \scrv$
and that $Z_{C_3}(q,v_1,v_2,v_3) > 0$ for $v_1,v_2,v_3 \in \scrv$.

The condition $\scrv \parallel \scrv \subseteq \scrv$ comes down to
\be
   (q/v_+) \parallel (q/v_+)  \;\le\;  v_+
\ee
or equivalently
\be
   v_+^3 - 2qv_+ - q^2  \;\ge\;  0   \;.
 \label{eq.diamond.inequality}
\ee
And since
$Z_{K_2^{(3)}}(q,v_1,v_2,v_3) = q [q + (v_1 \parallel v_2 \parallel v_3)]$,
the condition $Z_{K_2^{(3)}} < 0$ comes down to the two conditions
\begin{subeqnarray}
   v_+ \parallel v_+ \parallel v_+            & \le &   -q   \\
   (q/v_+) \parallel (q/v_+) \parallel v_+    & \le &   -q
\end{subeqnarray}
or equivalently
\begin{subeqnarray}
   (1+ v_+)^3            & \le &  1-q  \slabel{eq.cor3opt.a}   \\
   (1+ q/v_+)^2 (1+v_+)  & \le &  1-q  \slabel{eq.cor3opt.b}
\end{subeqnarray}
Expanding \reff{eq.cor3opt.b} leads to the condition that the cubic
\reff{eq.cor.0<q<1.block3optimal.cubic} must be $\le 0$,
so taking $v_+ = v_3^+(q)$ is permitted and in fact optimal.
Inequality \reff{eq.cor3opt.a} leads to the condition
$v_+ \le -1 + (1-q)^{1/3}$, which is weaker than  $v_+ \le v_3^+(q)$,
as can be seen by computing
\begin{eqnarray}
   v^3 + 3qv^2 + (q^2+2q)v + q^2 \bigg| _{v = -1 + (1-q)^{1/3}}
   & = &
   (1-q)^{4/3} \, \big[ 3 - 3(1-q)^{1/3} - q \big]    \nonumber \\[1mm]
   & > &
   (1-q)^{4/3} \, \big[ 3 - 3(1-q/3) - q \big]        \nonumber \\[1mm]
   & = &  
   0  \;.
\end{eqnarray}

Finally, let us show that $v_3^+(q)$ satisfies \reff{eq.diamond.inequality}.
Let $f(q,v)=v^3+3qv^2+(q^2+2q)v+q^2$ and $g(q,v)= v^3 - 2qv - q^2$.
Solving $f-g=0$ for $q$, we obtain $q=q_0(v)=-v(3v+4)(v+2)^{-1}$. 
Since $f(q_0(v),v)=v^4(v+1)(v+2)^{-2}\neq 0$ for all $v\in (-1,0)$, the 
curves $f=0$ and $g=0$ do not intersect when  $v\in (-1,0)$. It follows 
that $g(q,v_3^+(q))$ is non-zero with constant sign for all $q\in (0,1)$.
Taking $q=\alpha=(1+\sqrt{5})/4$, we have $v_3^+(\alpha)=-1/2$ and 
$g(\alpha,-1/2)=(\sqrt{5}-2)/8>0$. Thus \reff{eq.diamond.inequality} holds
when $v_+=v_3^+(q)$.
\qed

{\bf Remark.}
It is easy to show that
\be
    v_2^+(q) \;<\; v_3^+(q) \;<\; v_4^+(q)  \;.
\ee
Indeed, let us write $q=1-r^m$ with $0<r<1$
and substitute $v = v_m^+(q) = -(1-r^m)/(1+r)$ into the cubic
\reff{eq.cor.0<q<1.block3optimal.cubic};
we get $-r^2 (1-r)^3 < 0$ for $m=2$
and $r^3 (1-r)^4 (1+r^2)^2 > 0$ for $m=4$.
\qed

Here is the matroidal analogue of Corollary~\ref{cor.0<q<1.block34}:

\begin{corollary}
   \label{cor.0<q<1.block34.matroid}
Let $2 \le m \le 4$,
and let $M$ be a matroid with ground set $E$,
in which each 2-connected component contains at least $m$ edges.
Let $0 < q < 1$, and suppose that
\be
   -[1 + (1-q)^{1/m}]  \;<\;  v_e  \;<\;  - \, {q \over 1 + (1-q)^{1/m}}
 \label{eq.cor.0<q<1.block34.matroid}
\ee
for all $e \in E$.
Then $(-1)^{r(M)} \Ztilde_M(\qvbf) > 0$.
\end{corollary}

This is a nontrivial generalization of Corollary~\ref{cor.0<q<1.block34},
since 
there exists a non-graphic 2-connected matroid on four elements
that has to be included in the base case of the induction,
namely, the rank-2 uniform matroid $U_{2,4}$.

\proofof{Corollary~\ref{cor.0<q<1.block34.matroid}}
In addition to what was already done in proving
Corollary~\ref{cor.0<q<1.block34},
we need to prove that $\Ztilde_{U_{2,4}}(\qvbf) > 0$
whenever $0 < q < 1$ and $v_i \in I_4$ for $i=1,2,3,4$;
here $I_4$ is the self-dual interval defined in \reff{def.Im},
i.e.\ $I_4 = (v_-,v_+)$ where $v_- = -1-(1-q)^{1/4}$ and $v_+ = q/v_-$.
By \reff{eq.Imineq} we have $v_+ < -1 + (1-q)^{1/4}$
[this just says that $I_4 \subset I_4^{\rm cocyc}$],
hence $v_- + v_+ < -2$.

A simple computation gives
\be
   q^2 \Ztilde_{U_{2,4}}(\qvbf)
   \;=\;
   \prod_{1\leq i\leq 4}(1+v_i) \,+\,
     (q-1)\sum_{1\leq i\leq 4}v_i  \,+\,  (q^2-1)
   \;.
 \label{eq.Z.U24}
\ee
Since this is multiaffine and symmetric in the $\{v_i\}$,
it suffices to check that $q^2 \Ztilde_{U_{2,4}}(\qvbf) > 0$
for the five cases
$\bv = (v_-,v_-,v_-,v_-)$,
$(v_-,v_-,v_-,v_+)$,
$(v_-,v_-,v_+,v_+)$,
$(v_-,v_+,v_+,v_+)$ and
$(v_+,v_+,v_+,v_+)$.
Moreover, by self-duality of $U_{2,4}$ it suffices to check
the first three cases.
The first and third cases are handled by noting that
$\prod (1+v_i) > 0$ and $\sum v_i < -4$,
so that
\be
   q^2 \Ztilde_{U_{2,4}}(\qvbf)
   \;>\;
   (4-4q) + (q^2-1)
   \;=\;
   3 - 4q + q^2
   \;=\;
   (1-q)(3-q)
   \;>\;
   0   \;.
   \qquad
\ee
To handle the second case, let us define
$v'_+ \equiv -1 + (1-q)^{1/4} > v_+$ and compute
\begin{subeqnarray}
   q^2 \Ztilde_{U_{2,4}}(\qvbf)
   & = &
   (1+v_-)^3 (1+v_+) \,+\, (q-1)(3v_- + v_+) \,+\, (q^2-1)
      \qquad   \\[1mm]
   & > &
   (1+v_-)^3 (1+v'_+) \,+\, (q-1)(3v_- + v'_+) \,+\, (q^2-1)
   \;. \qquad
\end{subeqnarray}
Since $(1+v_-)^3 (1+v'_+) = -(1-q)$ and $3v_- + v'_+ < -4$, we have
\be
   q^2 \Ztilde_{U_{2,4}}(\qvbf)
   \;>\;
   -(1-q) + (4-4q) + (q^2-1)
   \;=\;
   2 - 3q + q^2
   \;=\;
   (1-q)(2-q)
   \;>\;
   0   \;.
\ee
\qed

Finally, Corollary~\ref{cor.0<q<1.block3optimal} extends immediately
to matroids:

\begin{corollary}
   \label{cor.0<q<1.block3optimal.matroid}
Let $M$ be a matroid with ground set $E$,
in which each 2-connected component contains at least three elements.
Let $0 < q < 1$, and suppose that
\be
   {q \over v_3^+(q)}  \;<\;  v_e  \;<\; v_3^+(q)
\ee
for all $e \in E$, where $v_3^+(q)$ is the unique real root
of the cubic equation
\be
   v^3 \,+\, 3qv^2 \,+\, (q^2+2q)v \,+\, q^2  \;=\;  0 \;.
 \label{eq.cor.0<q<1.block3optimal.cubic.matroid}
\ee
Then $(-1)^{r(M)} \Ztilde_M(\qvbf) > 0$.
\end{corollary}

\section{The diamond operation}\label{sec.diamond}

We will show in this section that the
hypotheses (a)--(c) of Theorem~\ref{thm.0<q<1.block.abstract}
imply that $q\leq 32/27$ and that $\scrv$ is contained in a particular interval
$I_\Diamond(q)$.
Our results will also be used in
Section~\ref{sec.q>1}
to obtain zero-free regions for $Z_G(q,v)$ when $q\in (1,32/27]$.
The ``diamond operation'',
in which one or more edges are replaced
by the parallel connection of two two-edge paths, will play a key role.
For any graph $G$, let us denote by $\Diamond(G)$
the graph in which {\em every}\/ edge of $G$ is replaced by a diamond.
And let us write
\be
   \Diamond_q(v)
   \;\equiv\;
   (v \series v) \parallel (v \series v)
   \;=\;
   {v^2 (v^2 + 4v + 2q)  \over  (q+2v)^2}
 \label{def.Diamond}
\ee
for the corresponding map of edge weights.
(This corresponds to the ``diagonal'' case,
 in which all four edges of the diamond get the {\em same}\/ weight $v$.)
Then, for any graph $G$ we have
\be
   Z_{\Diamond(G)}(q,v)
   \;=\;
   (q+2v)^{2|E(G)|} Z_G(q, \Diamond_q(v))
 \label{eq.diamond_identity}
\ee
when $v \neq -q/2$, by virtue of the series and parallel reduction rules
\reff{eq.series2}/\reff{eq.parallel2}.\footnote{
   These relations were found by various physicists in the early 1980s,
   in the course of work on Potts models on ``hierarchical lattices'':
   see e.g.\ \cite{Derrida_83,Itzykson_85}.
}
The case $v = -q/2$, which corresponds to $\Diamond_q(v) = +\infty$,
can be handled by a limiting process, using the fact that
$Z_G(q,w) \approx q^{k(G)} w^{|E(G)|}$ as $w \to \infty$;
combining this with \reff{def.Diamond}/\reff{eq.diamond_identity}, we obtain
\be
   Z_{\Diamond(G)}(q,-q/2)
   \;=\;
   q^{k(G)} (q/2)^{4|E(G)|}
   \;.
 \label{eq.diamond_identity.limiting}
\ee
In what follows we make the natural convention that
$\Diamond_q(+\infty) \equiv \lim\limits_{v \uparrow \infty} \Diamond_q(v)
 = +\infty$.

%
\begin{figure}[p]
  \vspace*{-5mm}
  \centering
  \includegraphics[width=3.8in]{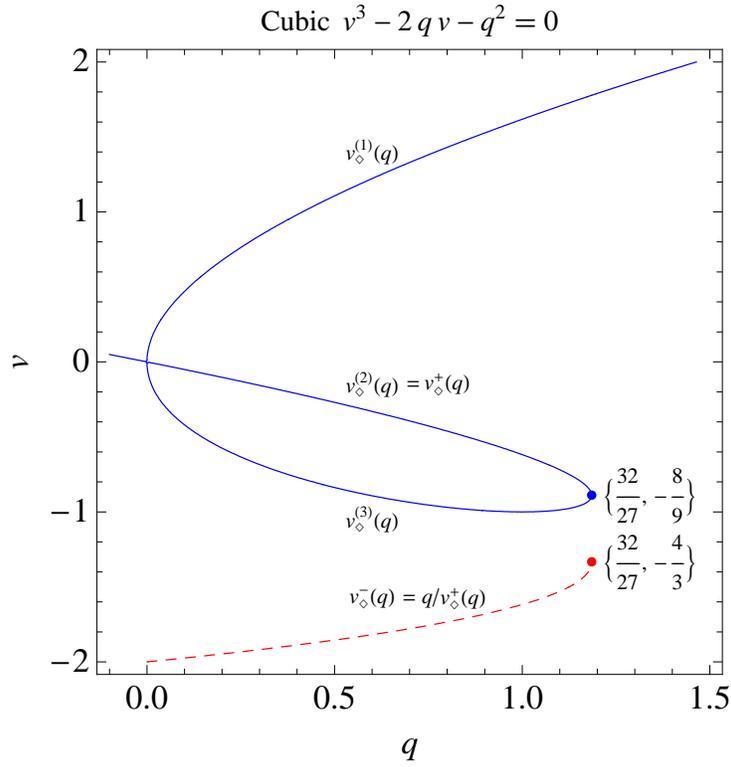}
  \caption{
      Solid blue curve shows the solutions of $v^3-2qv-q^2 = 0$.
      Dashed red curve shows the ``dual'' to the middle branch.
  }
\label{fig.diamondplot}
\end{figure}
%
%
\begin{figure}[p]
  \centering
  \includegraphics[width=4.2in]{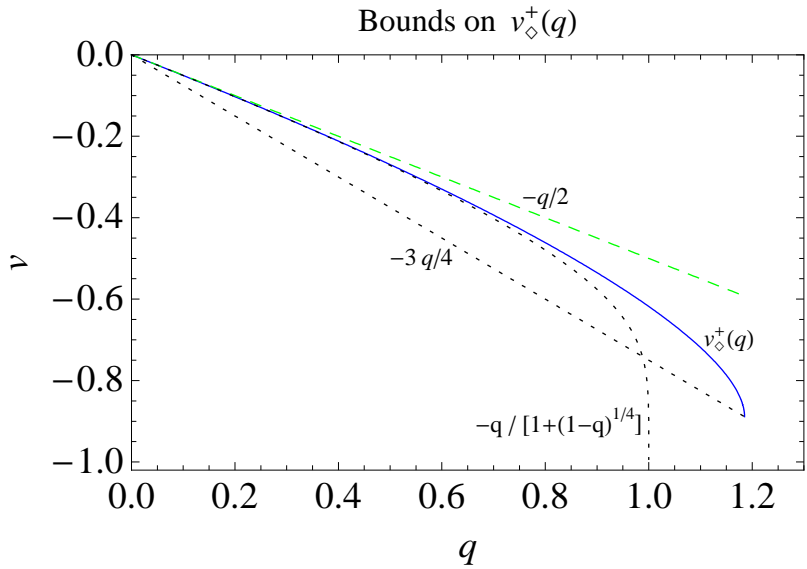}
  \caption{
      Solid blue curve shows $v_\Diamond^+(q)$.
      Dashed green line shows the upper bound $-q/2$.
      Dotted black curves show the lower bounds $-3q/4$
      and (for $0 \le q \le 1$) $-q/[1+(1-q)^{1/4}]$.
  }
\label{fig.diamondplot2}
\end{figure}
%
%
\begin{figure}[p]
  \centering
  \includegraphics[width=4.2in]{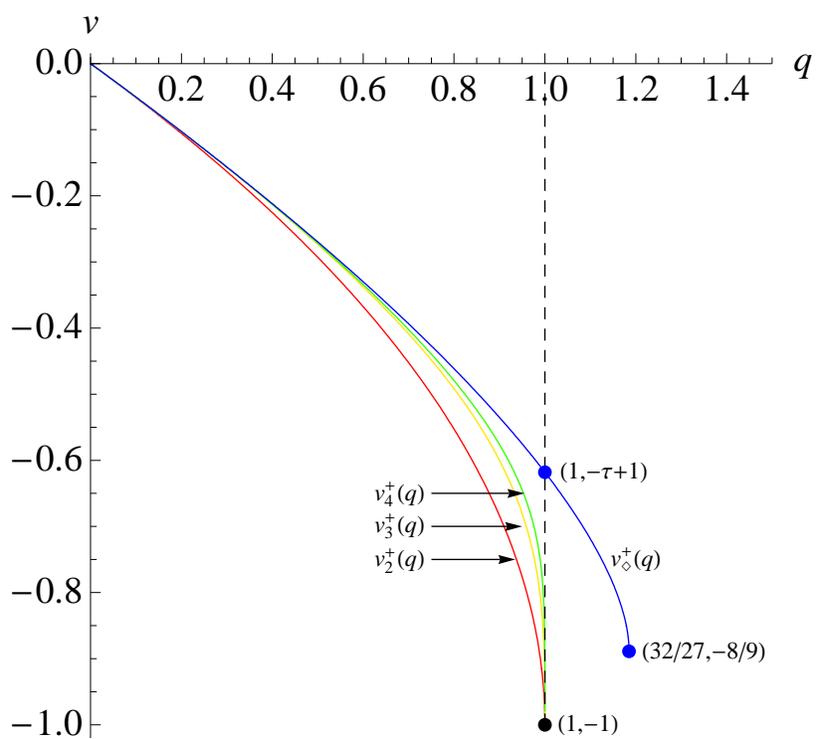}
  \caption{
      Solid blue curve shows $v_\Diamond^+(q)$.
      Solid red, yellow and green curves show
      $v_m^+(q)$ for $m=2,3,4$, respectively.
      Here $\tau = (1+\sqrt{5})/2$ is the golden ratio.
  }
\label{fig.plotdiaall}
\end{figure}
%

A central role in our analysis
will be played by the fixed points of the diamond map,
which satisfy $\Diamond_q(v) = v$ or equivalently
(excluding the trivial fixed points $v=0$ and $v=+\infty$)
\be
   v^3 - 2qv - q^2  \;=\;  0   \;.
 \label{eq.diamond.fixedpoint.cubic}
\ee
The cubic \reff{eq.diamond.fixedpoint.cubic} has one real root for $q < 0$,
three real roots for $0 < q < 32/27$, and one real root for $q > 32/27$
(see Figure~\ref{fig.diamondplot}).
For $0 < q \le 32/27$,
let us denote by $v_\Diamond^{(i)}$ ($i=1,2,3$)
the three roots of this cubic in decreasing order:
\be
   v_\Diamond^{(1)}  \:>\; 0  \;>\;  v_\Diamond^{(2)}
      \;\gtalmost\; v_\Diamond^{(3)}  \;\gtalmost\; -1
   \;,
\ee
where the first (resp.\ second) inequality $\gtalmost$ is strict
except at $q=32/27$ (resp.\ $q=1$).
We are especially interested in the middle branch $v_\Diamond^{(2)}$,
which we shall denote also by $v_\Diamond^+(q)$:
it decreases monotonically
from $v=0$ at $q=0$ to $v=-8/9$ at $q=32/27$,
and is given explicitly by the horrendous expression\footnote{
   We remark that the quantity
   $\smfrac{1}{16} [27q - 16 + i \sqrt{27q(32-27q)} \,]$
   lies for all $q \in [0,32/27]$ on the
   upper half of the unit circle in the complex plane;
   it runs from $-1$ at $q=0$ to $+1$ at $q=32/27$.
}
\begin{eqnarray}
  v_\Diamond^+(q)
   & = &
   \frac{3q}{2} \left[
      \{ \smfrac{1}{16} [27q - 16 + i \sqrt{27q(32-27q)} \,] \}^{1/3}
      \, e^{2\pi i/3}
   \;+\;
   \right.
       \nonumber \\
   & & \left.
   \quad \;+\;
     \{ \smfrac{1}{16} [27q - 16 + i \sqrt{27q(32-27q)}] \}^{-1/3}
     \,  e^{-2\pi i/3}
   \;-\;  1
   \right]^{-1}
   \qquad
\end{eqnarray}
or by the power series\footnote{
   This power series can be obtained by inserting
   $v=-(q/2)(1+w)$ into \reff{eq.diamond.fixedpoint.cubic}
   and using the Lagrange inversion formula to determine $w(q)$.
}
\be
   v_\Diamond^+(q)  \;=\;  -\frac{q}{2} \,-\,
       \sum_{m=1}^\infty \frac{1}{2m \, 8^m} {3m \choose m-1} \, q^{m+1}
   \;,
\ee
which is convergent for $|q| < 32/27$
and shows that all derivatives of $v_\Diamond^+(q)$
are strictly negative for $0 < q < 32/27$.
Putting $f(q,v)=v^3 - 2qv - q^2$, we have $f(q,-3q/4)\geq 0 >f(q,-q/2)$
for all $0<q\leq 32/27$ and hence
\begin{equation}
    -3q/4 \;\le\; v_\Diamond^+(q) \;<\; -q/2
 \label{concave}
\end{equation}
[these bounds alternatively follow from the concavity of $v_\Diamond^+(q)$].
For $0 < q < 1$ we also have\footnote{
   Writing $q = 1-r^4$ with $0 < r < 1$,
   we find $f(q, -q/[1 + (1-q)^{1/4}]) = r(1-r)^4 (1+r^2)^2 > 0$
   for $0 < q < 1$.
   At $q=1$ this vanishes because $-q/[1 + (1-q)^{1/4}]$
   touches the {\em bottom}\/ branch $v_\Diamond^{(3)}(q)$.
}
$f(q, -q/[1 + (1-q)^{1/4}]) > 0$
and hence
\be
   v_4^+(q)  \;\equiv\; -q/[1 + (1-q)^{1/4}]  \;<\; v_\Diamond^+(q)
   \;.
\ee 
These bounds are illustrated in Figure~\ref{fig.diamondplot2}.
In Figure~\ref{fig.plotdiaall} we compare $v_\Diamond^+(q)$
with the functions $v_2^+(q)$, $v_3^+(q)$ and $v_4^+(q)$
introduced in the preceding section.

The fixed point $v_\Diamond^+(q)$ is repulsive for $0 < q < 32/27$
and becomes marginal at $q = 32/27$:
more precisely, the ``multiplier''
\be
   \lambda_\Diamond(q)
   \;\equiv\;
   \left. {d \Diamond_q(v) \over dv} \right| _{v = v_\Diamond^+(q)}
\ee
decreases monotonically from $\lambda_\Diamond = +\infty$ at $q=0$
to $\lambda_\Diamond = +1$ at $q=32/27$.

We then define $v_\Diamond^-(q)$ to be the dual point
\begin{subeqnarray}
   v_\Diamond^-(q)
   & \equiv &  {q \over v_\Diamond^+(q)}
        \\[3mm]
   & = &
   \frac{2}{3} \left[
      \{ \smfrac{1}{16} [27q - 16 + i \sqrt{27q(32-27q)} \,] \}^{1/3}
      \, e^{2\pi i/3}
   \;+\;
   \right.
       \nonumber \\
   & & \left.
   \quad +\;
     \{ \smfrac{1}{16} [27q - 16 + i \sqrt{27q(32-27q)}] \}^{-1/3}
     \,  e^{-2\pi i/3}
   \;-\;  1
   \right]
   \qquad  \\[2mm]
   & = &
   -2 \,+\,
       \sum_{m=1}^\infty \frac{2}{m \, 8^m} {3m-2 \choose m-1} \, q^{m}
     \;,
\end{subeqnarray}
which increases monotonically, with all derivatives nonnegative,
from $v=-2$ at $q=0$ to $v=-4/3$ at $q=32/27$
(see again Figure~\ref{fig.diamondplot}).\footnote{
   The power series for $v_\Diamond^-(q)$ can be obtained by inserting
   $v=-(q/2)(1+w)$ into \reff{eq.diamond.fixedpoint.cubic}
   and using the Lagrange inversion formula to determine
   $q/v(q) = -2/[1+w(q)]$.
   This series is manifestly convergent for $|q| < 32/27$,
   and every derivative of $v_\Diamond^-(q)$ is
   strictly positive for $0 < q < 32/27$.
}
Finally, we let $I_\Diamond(q)$ be the ``diamond interval''
\be
   I_\Diamond(q)   \;=\;  \big[ v_\Diamond^-(q), \, v_\Diamond^+(q)  \big]
   \;.
\ee

The key facts about $v_\Diamond^\pm(q)$ are summarized in the following lemma,
which will be proven at the end of this section:

\begin{lemma}
  \label{lemma.diamond1}
Let $0 < q \le 32/27$.  Then:
\begin{itemize}
   \item[(a)]  $v_\Diamond^+(q) \series v_\Diamond^+(q) = v_\Diamond^-(q)$.
   \item[(b)]  $v_\Diamond^-(q) \parallel v_\Diamond^-(q) = v_\Diamond^+(q)$.
   \item[(c)]  $I_\Diamond(q)$ is self-dual, i.e.\ it is invariant under
       $v \mapsto q/v$.
   \item[(d)]  $I_\Diamond(q) \parallel I_\Diamond(q) \subseteq I_\Diamond(q)$.
   \item[(e)]  $I_\Diamond(q) \series I_\Diamond(q) \subseteq I_\Diamond(q)$.
\end{itemize}
\end{lemma}

As a strong converse to Lemma~\ref{lemma.diamond1}(d,e),
we have the following {\em necessary}\/ condition for
invariance under parallel and series connection:

\begin{proposition}
   \label{prop.diamond1.converse}
Let $q > 0$ and let $\emptyset \neq \scrv \subset \R$ satisfy
\begin{itemize}
   \item[(a)]  $\scrv \subseteq (-\infty,0)$
   \item[(b)]  If $v \in \scrv$, then $v \parallel v \in \scrv$.
   \item[(c)]  If $v \in \scrv$, then $v \series v \in \scrv$.
\end{itemize}
Then $q \le 32/27$ and $\scrv \subseteq I_\Diamond(q)$.
\end{proposition}

\noindent
Please note that hypotheses (b) and (c) are weaker than
$\scrv \parallel \scrv \subseteq \scrv$ and
$\scrv \series \scrv \subseteq \scrv$,
as they require invariance only under ``diagonal''
parallel and series connection.  Thus Proposition
\ref{prop.diamond1.converse} immediately implies:

\begin{corollary}\label{propdiamond}
Suppose that $q > 0$ and $\scrv \neq \emptyset$ satisfies hypotheses (a)--(c)
of Theorem~\ref{thm.0<q<1.block.abstract}.
Then $q\leq 32/27$ and $\scrv \subseteq I_\Diamond(q)$.
\end{corollary}

In the special case of {\em self-dual}\/ intervals
$\scrv = (q/v_+, v_+)$ with $v_+ < 0$,
we can give a {\em necessary and sufficient}\/ condition
to have invariance under parallel or series connection:

\begin{proposition}
   \label{prop.diamond1.converse2}
Let $q > 0$ and let $\scrv = (q/v_+, v_+)$ with $-\sqrt{q} < v_+ < 0$
(so that $\scrv \neq \emptyset$).
Then the following are equivalent:
\begin{itemize}
   \item[(a)]  $\scrv \parallel \scrv \subseteq \scrv$
   \item[(b)]  $\scrv \series \scrv \subseteq \scrv$
   \item[(c)]  $q \le 32/27$ and $\max[v_\Diamond^{(3)}(q),\, -q] \le v_+ \le
                 v_\Diamond^{(2)}(q)$.
\end{itemize}
Furthermore, we have
\be
   \max[v_\Diamond^{(3)}(q),\, -q]
   \;=\;
   \cases{-q                     & \hbox{\rm if } $0 < q \le 1$  \cr
          \noalign{\vskip 3mm}
          v_\Diamond^{(3)}(q)    & \hbox{\rm if } $1 \le q \le 32/27$ \cr
         }
  \label{eq.prop.diamond1.converse2}
\ee
\end{proposition}

The following further facts are relevant to the applicability
of Theorem~\ref{thm.0<q<1.block.abstract}:

\begin{lemma}
  \label{lemma.diamond2}
\hfill\break
\vspace*{-0.6cm}
\par\noindent
\begin{itemize}
   \item[(a)]  If $q > 0$ and $v > 0$, then $\Diamond_q^k(v) > 0$ for all $k$.
   \item[(b)]  If $0 < q \le 32/27$ and $v_\Diamond^+(q) < v < 0$,
      then $\Diamond_q(v) > v$,
      and $\Diamond_q^k(v) \ge 0$ for all sufficiently large $k$.
   \item[(c)]  If $q > 32/27$ and $v < 0$,
      then $\Diamond_q(v) > v$,
      and $\Diamond_q^k(v) \ge 0$ for all sufficiently large $k$.
   \item[(d)]  If $q \ge 1$ and $v \in \R$,
      then $\Diamond_q(v) \ge -q$,
      with strict inequality except when $(q,v)=(1,-1)$.
\end{itemize}
\end{lemma}

\noindent
It follows that in cases (b) and (c),
the sequence $\{ \Diamond_q^k(v) \}_{k \ge 0}$
is strictly increasing as long as it stays negative;
and once it goes nonnegative, it stays nonnegative
(but need no longer be increasing).
Note also that, in cases (b)--(d),
if one iterate $\Diamond_q^k(v)$ ($k \ge 0$) happens to equal $-q/2$,
then the next iterate and all subsequent iterates will equal $+\infty$
(which is indeed $\ge 0$).

\begin{corollary}
   \label{cor.diamond2}
Let $G$ be a graph.
\begin{itemize}
   \item[(a)] If $0<q\leq 32/27$ and $v>v_\Diamond^+(q)$, then 
       $Z_{\Diamond^k(G)}(q,v) > 0$ for all sufficiently large $k$.
   \item[(b)] If $0<q\leq 32/27$ and $v<v_\Diamond^-(q)$, then 
       $Z_{\Diamond^k(G)^{(2)}}(q,v) > 0$ for all sufficiently large $k$.
[Here $G^{(2)}$ denotes the graph obtained from $G$
 by replacing each edge by two parallel edges.]
   \item[(c)] If $q> 32/27$ and $v\in \R$, then 
       $Z_{\Diamond^k(G)}(q,v) > 0$ for all sufficiently large $k$.
\end{itemize}
\end{corollary}

We have already seen in Corollary~\ref{propdiamond}
that if $q > 0$ and $\scrv \neq \emptyset$ satisfy hypotheses (a)--(c)
of Theorem~\ref{thm.0<q<1.block.abstract},
then $q\leq 32/27$ and $\scrv \subseteq I_\Diamond(q)$.
We can now show, using Corollary~\ref{cor.diamond2},
that if $m \ge 2$, $\gamma \in \{0,1\}$, $q > 0$ and $\scrv \neq \emptyset$
satisfy the {\em conclusion}\/ of Theorem~\ref{thm.0<q<1.block.abstract}
(with $\scrg \supseteq$ series-parallel graphs), 
then we must either have
\begin{itemize}
   \item[(a)] $\gamma=0$, $q < 1$ and $\scrv \subseteq I_\Diamond(q)$
\end{itemize}
or else
\begin{itemize}
   \item[(b)] $\gamma=1$, $1 < q \le 32/27$ and
      $\scrv \subseteq I_\Diamond(q)$
\end{itemize}
--- and this is so no matter how large we take $m$ to be.

\begin{corollary}
   \label{cor.diamond2a}
Fix $m \ge 2$, $\gamma\in\{0,1\}$ and $q > 0$.
Suppose that $\scrv \subset \R$ ($\scrv \neq \emptyset$)
satisfies the conclusion of Theorem~\ref{thm.0<q<1.block.abstract}
for some class $\scrg \supseteq$ series-parallel graphs.  Then either:
\begin{itemize}
   \item[(a)] $\gamma=0$, $q < 1$ and $\scrv \subseteq I_\Diamond(q)$; or
   \item[(b)] $\gamma=1$, $1 < q \le 32/27$ and
         $\scrv \subseteq I_\Diamond(q)$.
\end{itemize}
\end{corollary}

%

Let us now prove all these results:

\proofof{Lemma~\ref{lemma.diamond1}}
(a) The equation $v \series v = q/v$ is precisely the
cubic equation \reff{eq.diamond.fixedpoint.cubic}
satisfied by $v = v_\Diamond^+(q)$.

(b) follows from (a) by duality.

(c) is easy.

(d) Since $-2 < v_\Diamond^-(q) < -1 < v_\Diamond^+(q) < 0$ by
(\ref{concave}),
the inequalities
$-1 < v_\Diamond^+(q) \parallel v_\Diamond^+(q) < v_\Diamond^+(q)$
and
$v_\Diamond^-(q) < v_\Diamond^+(q) \parallel v_\Diamond^-(q) < -1$
are trivial.
And by (b) we have the equality
$v_\Diamond^-(q) \parallel v_\Diamond^-(q) = v_\Diamond^+(q)$.

(e) follows from (d) by duality.
\qed

\proofof{Lemma~\ref{lemma.diamond2}}
(a) is obvious.

(b,c) If either $0 < q \le 32/27$ and $v_\Diamond^+(q) < v < 0$,
or $q > 32/27$ and $v < 0$,
then the cubic \reff{eq.diamond.fixedpoint.cubic}
has the sign $v^3 - 2qv - q^2 < 0$.
Since
\be
   \Diamond_q(v) - v
   \;=\;
   {v (v^3 - 2qv - q^2)  \over (q+2v)^2}
   \;,
\ee
it follows that $\Diamond_q(v) > v$ in these cases.
(The value is unambiguously $+\infty$ if $v = -q/2$.)

Let us next prove that $\Diamond_q^k(v) \ge 0$ for all sufficiently large $k$.
Note first that if $\Diamond_q^k(v) \ge 0$,
then $\Diamond_q^\ell(v) \ge 0$ for all $\ell \ge k$,
since $\Diamond_q$ obviously maps $[0,+\infty]$ into itself.
So it suffices to prove that $\Diamond_q^k(v) \ge 0$ for at least one $k$.
Assume the contrary:
then, by virtue of what has already been shown, we have
\be
   v \,<\, \Diamond_q(v) \,<\, \Diamond_q^2(v) \,<\,
           \Diamond_q^3(v) \,<\, \ldots \,<\, 0
\ee
and hence the sequence $\Diamond_q^k(v)$ tends to a limit $v_*$
satisfying $v < v_* \le 0$
(and in particular $v_* > v_\Diamond^+(q)$ in case $q \le 32/27$).
Since $v_*$ must be a fixed point of $\Diamond_q$,
the only possibility is $v_* = 0$.
But for $-q/2 < v < 0$ we have $\Diamond_q(v) > 0$,
which rules out the possibility that $\Diamond_q^k(v)$ tends to 0 from below.

(d) follows immediately from
\be
   \Diamond_q(v) + q
   \;=\;
   { (1+v)^4 \,+\, (q-1) \Big[ 6 \big( v + {q+1 \over 3} \big)^2 +
                               {1 \over 3} (q^2 - q + 1)
                         \Big]
    \over
    (q+2v)^2
   }
   \;.
\ee
(The value is unambiguously $+\infty$ if $v = -q/2$.)
\qed

\proofof{Corollary~\ref{cor.diamond2}}
Since we are asserting that $Z_G > 0$,
it suffices to consider connected graphs $G$.

(a) We have $\Diamond_q^k(v) \ge 0$ for all sufficiently large $k$
by Lemma~\ref{lemma.diamond2}(a,b).
So consider such a $k$.
If none of the iterates
$v, \Diamond_q(v), \Diamond_q^2(v), \ldots, \Diamond_q^{k-1}(v)$
happens to equal $-q/2$,
then it follows from \reff{eq.diamond_identity} that, for any graph $G$,
the {\em bivariate}\/ Tutte polynomial
of the graph $\Diamond^k(G)$ satisfies
\be
   Z_{\Diamond^k(G)}(q,v)  \;=\;
   \hbox{positive prefactors} \,\times\, Z_G(q, \Diamond_q^k(v))
   \;>\; 0
   \;.
\ee 
If, on the other hand, one of the iterates $\Diamond_q^\ell(v)$
[with $0 \le \ell \le k-1$]
equals $-q/2$, then it follows from
\reff{eq.diamond_identity} and \reff{eq.diamond_identity.limiting} that
\be
   Z_{\Diamond^k(G)}(q,v)  \;=\;
   \hbox{positive prefactors} \,\times\,
      q \, (q/2)^{4 |E(\Diamond^{k-\ell-1}(G))|}
   \;>\; 0
   \;.
\ee

(b) is an immediate consequence of (a),
since $v < v_\Diamond^-(q)$ implies
$v \parallel v > v_\Diamond^+(q)$.

(c) For all $q > 32/27$ and $v \in \R$, we have
$\Diamond_q^k(v) \ge 0$ for all sufficiently large $k$
by Lemma~\ref{lemma.diamond2}(a,c);
so $Z_{\Diamond^k(G)}(q,v) > 0$
by the same argument as in part (a).
\qed

\proofof{Corollary~\ref{cor.diamond2a}}
We may show that either $\gamma = 0$ and $q < 1$
or else $\gamma = 1$ and $q > 1$ by an argument
similar to that used in the proof of Proposition~\ref{gamma}
(we leave the details to the reader).
To show that $q \le 32/27$ and $\scrv \subseteq I_\Diamond(q)$,
suppose the contrary:
then we use Corollary~\ref{cor.diamond2} with $G=K_2$ and $K_3$
to construct 2-connected series-parallel graphs $H$,
with an arbitrarily large number of edges,
whose vertex-set sizes have {\em both}\/ parities.
(Here we have used the fact that if $G$ is a non-separable graph,
 then $\Diamond^k(G)$ and $\Diamond^k(G)^{(2)}$ are both non-separable,
 and the parity of the size of their vertex sets
 is the same as that of $G$.)
Since Corollary~\ref{cor.diamond2} yields $Z_H > 0$
while the conclusion of Theorem~\ref{thm.0<q<1.block.abstract}
asserts that $(-1)^{|V(H)|-1+\gamma} Z_H > 0$,
one of the two parities yields a counterexample.
\qed

\proofof{Proposition~\ref{prop.diamond1.converse}}
Hypotheses (a)--(c) imply in particular
that $\Diamond_q(\scrv) \subseteq \scrv \subseteq (-\infty,0)$.
By Lemma~\ref{lemma.diamond2}(a,b),
this is possible with $\scrv \neq \emptyset$ only if $q \le 32/27$
and $\scrv \subseteq (-\infty, v_\Diamond^+(q)]$.
On the other hand, if $v < v_\Diamond^-(q)$,
then [since $v_\Diamond^-(q) \le -1$] we have
$v \parallel v > v_\Diamond^-(q) \parallel v_\Diamond^-(q) = v_\Diamond^+(q)$.
So, by hypothesis (b),
we must have $\scrv \subseteq [v_\Diamond^-(q), v_\Diamond^+(q)]$.
\qed

\proofof{Proposition~\ref{prop.diamond1.converse2}}
Since $\scrv = (q/v_+, v_+)$ is self-dual,
$\scrv \parallel \scrv \subseteq \scrv$
is equivalent to $\scrv \series \scrv \subseteq \scrv$;
so let us check the former.
An obvious necessary condition is $q/v_+ \le -1 \le v_+$,
i.e.\ $v_+ \ge \max(-1,-q)$.
If these conditions are satisfied, a necessary and sufficient condition
is then
\be
   (q/v_+) \parallel (q/v_+)  \;\le\;  v_+
\ee
or equivalently
\be
   v_+^3 - 2qv_+ - q^2  \;\ge\;  0   \;.
\ee
But this is just the ``diamond cubic'', so we must have
either $v_+ \ge v_\Diamond^{(1)}(q) > 0$
or else $v_\Diamond^{(2)}(q) \ge v_+ \ge v_\Diamond^{(3)}(q) \ge -1$.

Finally, let us prove \reff{eq.prop.diamond1.converse2}.
We have $v^3 - 2qv - q^2 \big| _{v=-q} = q^2 - q^3$.
For $0<q<1$ this is positive,
so we must have $v_\Diamond^{(3)}(q) < -q < v_\Diamond^{(2)}(q)$.
For $1 < q \le 32/27$ this is negative,
so we must have either $-q < v_\Diamond^{(3)}(q)$
or $v_\Diamond^{(2)}(q) < -q < v_\Diamond^{(1)}(q)$;
but the latter is excluded because we know that
$-3q/4 \le v_\Diamond^{(2)}(q)$ by (\ref{concave}).
\qed



\section{The interval \boldmath $q\in (1, 32/27]$}   \label{sec.q>1}

Let $G$ be a loopless graph with $n$ vertices, $c$ components,
and $b$ nontrivial blocks\footnote{
   Let us recall that we call a block {\em trivial}\/
   if it has only one vertex, and {\em nontrivial}\/ otherwise.
};
then Theorem~\ref{thm1.1}(e) states that
$P_G(q)$ is nonzero with sign $(-1)^{n+c+b}$ for $1 < q \le 32/27$
\cite{Jackson_93}.
An analogous result also holds for loopless matroids \cite{Edwards_98}.
In this section we shall use Theorem ~\ref{thm.0<q<1.block.abstract}
to generalize these results
to the multivariate Tutte polynomial.
The multivariate approach allows us to replace the detailed
graph-theoretic proof of
\cite{Jackson_93} by a much simpler proof involving elementary
calculus.

We need to find intervals
$\scrv = (v_-, v_+)$ or $\scrv = [v_-, v_+]$
satisfying hypotheses (a)--(d) of Theorem ~\ref{thm.0<q<1.block.abstract}
with $\gamma=1$,
for the case $\scrg = $ all graphs.
For simplicity, let us restrict attention to {\em self-dual}\/ intervals,
i.e.\ $v_- = q/v_+$.
Then invariance under parallel connection
is equivalent to invariance under series connection;
and by Proposition~\ref{prop.diamond1.converse2},
these properties hold if and only if $q \le 32/27$
and
\be
    \left\{\!\! \begin{array}{ll}
                   -q                   & \hbox{if } 0 < q \le 1   \\[1mm]
                   v_\Diamond^{(3)}(q)  & \hbox{if } 1 \le q \le 32/27
                \end{array}
    \!\!\right\}
    \;\le\;  v_+  \;\le\;  v_\Diamond^{(2)}(q)
    \;.
\ee
Since $v_\Diamond^{(2)}(q) < -q/2$ by (\ref{concave}),
hypothesis (a) is then satisfied as well.
Finally, by Proposition~\ref{gamma},
we can restrict attention to the case $q>1$.

It remains to determine the conditions under which
also hypothesis (d) holds.
We have been able to do this,
and thus to find the optimal self-dual interval,
for the cases $m=2$ and $m=3$.

\bigskip

{\bf Case \boldmath $m=2$.}
The only non-separable graph with two edges is $K_2^{(2)} = C_2$.
We want to have $Z_{K_2^{(2)}}(q,v_1,v_2) = q [q + (v_1 \parallel v_2)] > 0$
for all $v_1,v_2 \in \scrv$,
hence we need $\scrv \parallel \scrv \subseteq (-q,\infty)$
[actually it will be $\subseteq (-q,0)$].
Since $q > 1$, the conditions $v_+ \parallel v_+ > -1 > -q$
and $(q/v_+) \parallel (q/v_+) -1 > -q$ hold trivially.
So the only nontrivial condition is
\be
   v_+ \parallel (q/v_+)  \;\ge\;  -q
\ee
or equivalently
\be
   v_+^2 + 2qv_+ + q  \;\le\;  0   \;.
 \label{eq.quadratic.m=2}
\ee
This means that $v_+$ must lie between the two roots of the quadratic
\reff{eq.quadratic.m=2}, which are $-q \pm \sqrt{q^2 - q}$.
Of course, we must also make sure that
$v_\Diamond^{(3)}(q) \le v_+ \le v_\Diamond^{(2)}(q)$
to satisfy hypotheses (a)--(c).
The maximal choice $v_+ = -q + \sqrt{q^2 - q}$
works whenever $1 < q \le 9/8$.\footnote{
   {\sc Proof:}  It suffices to check that
   $$ v^3 - 2qv - q^2 \bigg| _{v = -q + \sqrt{q^2 - q}}
      \;\ge\; 0  \;.$$
   But this equals $q \sqrt{q^2-q} \: [ 4q-3 - 4 \sqrt{q^2-q}]$,
   so we need $4q-3 \ge 4 \sqrt{q^2-q}$,
   i.e.\ $q \le 9/8$.
}
Otherwise the best we can do is to take
$v_+ = v_\Diamond^{(2)}(q) \equiv v_\Diamond^+(q)$.
We have therefore proven:

\begin{corollary}
   \label{cor.q>1.block.m=2}
Let $1 < q \le 32/27$ and define
\be
   \scrv_2  \;=\;
   \cases{  \big( -q - \sqrt{q^2 - q}, \: -q + \sqrt{q^2 - q} \, \big)
               & if $1 < q \le 9/8$  \cr
            \noalign{\vskip 3mm}
            I_\Diamond(q) \,\equiv\, [v_\Diamond^-(q), v_\Diamond^+(q)]
               & if $9/8 \le q \le 32/27$  \cr
         }
\ee
Then $(-1)^{n+c+b} Z_G(\qvbf) > 0$
whenever $G=(V,E)$ is a loopless bridgeless graph
with $n$ vertices, $c$ components and $b$ nontrivial blocks,
and $v_e \in \scrv_2$ for all $e \in E$.
\end{corollary}

The interval $\scrv_2$ is the best possible {\em self-dual}\/
interval for Corollary~\ref{cor.q>1.block.m=2}, in the following senses:

(i) For all $q > 1$, the graph $K_2^{(2)}$ has a {\em multivariate}\/ root
at $v_1 = -q - \sqrt{q^2 - q}$, $v_2 = -q + \sqrt{q^2 - q}$.

(ii) For $1 < q \le 32/27$, $v > v_\Diamond^+(q)$, 
and $G$ an arbitrary  graph,
the graph $\Diamond^k(G)$ satisfies $Z_{\Diamond^k(G)}(q,v)>0$ by
Corollary \ref{cor.diamond2}(a).
This has the {\em wrong}\/ sign for
Corollary~\ref{cor.q>1.block.m=2} when $G$ is
 2-connected with an {\em odd}\/ number of vertices.
A similar argument holds if $1 < q \le 32/27$ and $v < v_\Diamond^-(q)$,
or $q>32/27$ and $v\in \R$, using Corollary \ref{cor.diamond2}(b) and (c),
respectively.
%
%


\bigskip

{\bf Case \boldmath $m=3$.}
The only non-separable graphs with two edges are $K_2^{(3)}$
and its dual $C_3$.
By self-duality, it suffices to consider the former.
We want to have
$Z_{K_2^{(3)}}(q,v_1,v_2,v_3) = q [q + (v_1 \parallel v_2 \parallel v_3)] > 0$
for all $v_1,v_2,v_3 \in \scrv$,
hence we need $\scrv \parallel \scrv \parallel \scrv \subseteq (-q,\infty)$
[again it will actually be $\subseteq (-q,0)$].
Since $q > 1$, it is easy to see that $|1+q/v_+| > |1+v_+|$,
i.e.\ the interval $\scrv$ extends farther to the left of $-1$
than to the right.
Therefore, the necessary and sufficient condition to have
$\scrv \parallel \scrv \parallel \scrv \subseteq (-q,\infty)$ is simply
\be
   (q/v_+) \parallel (q/v_+) \parallel (q/v_+) \;\ge\;  -q
\ee
or equivalently
\be
   v_+^3 + 3v_+^2 + 3qv_+ + q^2  \;\le\;  0   \;.
 \label{eq.cubic.m=3}
\ee
For $q > 1$ this cubic has a single real root
$v_+ = -1 + (q-1)^{1/3} - (q-1)^{2/3}$,\footnote{
   The derivative of this cubic,
   namely $3v^2 + 6v + 3q$, has discriminant $36-36q < 0$,
   so the cubic has strictly positive derivative on all of $\R$.
   It is easily verified by substitution that
   $v = -1 + (q-1)^{1/3} - (q-1)^{2/3}$ is indeed the root.
}
so the inequality \reff{eq.cubic.m=3} reduces to
\be
   v_+  \;\le\;  -1 \,+\, (q-1)^{1/3} \,-\, (q-1)^{2/3}
   \;.
\ee
Of course, we must also make sure that
$v_\Diamond^{(3)}(q) \le v_+ \le v_\Diamond^{(2)}(q)$
in order to satisfy hypotheses (a)--(c).
The maximal choice $v_+ = -1 + (q-1)^{1/3} - (q-1)^{2/3}$
works whenever $1 < q \le 9/8$.\footnote{
   {\sc Proof:}  It suffices to check that
   $$ v^3 - 2qv - q^2 \bigg| _{v = -1 + (q-1)^{1/3} - (q-1)^{2/3}}
      \;\ge\; 0  \;.$$
   Making the change of variables $r = (q-1)^{1/3}$,
   a short calculation shows that we need
   $r (1-2r) (r^2 - r +1)^2 \ge 0$,
   i.e.\ $r \le 1/2$, hence $q \le 9/8$.
}
(It is an amazing coincidence ---
 for which we have no deep explanation ---
 that both $m=2$ and $m=3$ give rise to the same
 crossover point $q=9/8$.)
Otherwise the best we can do is to take
$v_+ = v_\Diamond^{(2)}(q) \equiv v_\Diamond^+(q)$.
We have therefore proven:

\begin{corollary}
   \label{cor.q>1.block.m=3}
Let $1 < q \le 32/27$ and define
\be
   \scrv_3  \;=\;
   \cases{  \biggl( \, \displaystyle {q \over -1 + (q-1)^{1/3} - (q-1)^{2/3}}
                    , \,
                    -1 + (q-1)^{1/3} - (q-1)^{2/3} 
            \biggr)
               & if $1 < q \le 9/8$  \cr
            \noalign{\vskip 2mm}
            I_\Diamond(q) \,\equiv\, [v_\Diamond^-(q), v_\Diamond^+(q)]
               & if $9/8 \le q \le 32/27$  \cr
            \noalign{\vskip 1pt}
         }
\ee
Then $(-1)^{n+c+b} Z_G(\qvbf) > 0$
whenever $G=(V,E)$ is a
graph with $n$ vertices, $c$ components and $b$ blocks,
in which each block contains at least three edges,
and $v_e \in \scrv_3$ for all $e \in E$.
\end{corollary}

We can show that the interval $\scrv_3$ is best possible in the same way 
as for the
interval $\scrv_2$ of Corollary \ref{cor.q>1.block.m=2}.
Suppose $1 < q \le 9/8$. If either 
$v_1 = v_2 = v_3 < q/[-1 + (q-1)^{1/3} - (q-1)^{2/3}]$ 
and $G=K_2^{(3)}$, or $v_1 = v_2 = v_3 < q/[-1 + (q-1)^{1/3} - (q-1)^{2/3}]$,
and $G=C_3$, then $Z_G(q,\qvbf)$ has the wrong sign.
So the interval $\scrv_3$ is best possible for $1 < q \le 9/8$,
even in the univariate case, and even allowing subsets $\scrv$
that are not necessarily intervals and not necessarily self-dual.
When $q>9/8$, the argument is the same as that given for $\scrv_2$.

\bigskip

We close this section by extending our results to matroids.
Since all matroids with at most three elements are graphic,
Corollaries~\ref{cor.q>1.block.m=2} and \ref{cor.q>1.block.m=3}
and Theorem \ref{thm.0<q<1.block.abstract.matroid}
immediately imply:

\begin{corollary}
Let $1 < q \le 32/27$ and $m\in \{2,3\}$,
and let ${\cal V}_2$ and ${\cal V}_3$ be the intervals defined in
Corollaries~\ref{cor.q>1.block.m=2} and \ref{cor.q>1.block.m=3}, respectively.
Then $(-1)^{r(M)+b} \Ztilde_M(\qvbf) > 0$
whenever $M$ is a matroid of rank $r(M)$ on the ground set $E$
that has $b$ $2$-connected components,
in which each $2$-connected component contains at least $m$ elements,
and $v_e \in \scrv_m$ for all $e \in E$.
\end{corollary}

%
%
%
%
%
%
%

\section{Further refinements?}  \label{sec.further}

Let us conclude by making some remarks on the possibility
of extending the zero-free regions obtained in Sections~\ref{sec4.2} and
\ref{sec.q>1} to larger values of $m$.
We {\em conjecture}\/ that for $0<q<1$ and $1<q\leq 32/27$,
there exists a strictly increasing
family of self-dual intervals
$\scrv_m(q) = (q/v_m^+(q), v_m^+(q))$
that satisfy the hypotheses of Theorem~\ref{thm.0<q<1.block.abstract}
for the case $\scrg = $ all graphs,
and such that $\lim\limits_{m \to\infty} \scrv_m(q) =I_\Diamond(q)$.
If true, this would imply:

%
\begin{figure}[t]
  \vspace*{-5mm}
  \centering
  \includegraphics[width=5.5in]{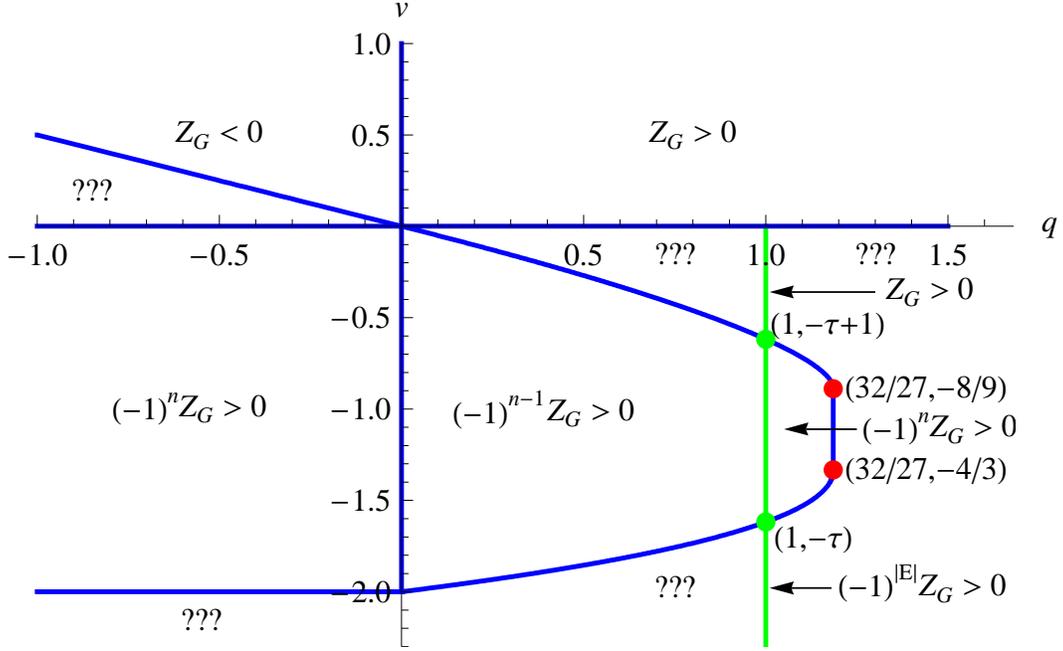}
  \caption{
      Conjectured limiting regions as $m \to\infty$
      in which the sign of $Z_G(q,{\bf v})$ can be controlled
      for loopless 2-connected graphs $G$ with at least $m$ edges.
      Here $\tau = (1+\sqrt{5})/2$ is the golden ratio.
  }
\label{fig.finalplot}
\end{figure}
%

\begin{conjecture}\label{conq<1}
Suppose $0<q<1$. Then there exists
a strictly increasing sequence of self-dual intervals $\scrv_m(q)$,
$m\geq 2$, such that
\begin{itemize}
   \item[(a)] $\lim\limits_{m\to \infty} \scrv_m(q)=I_\Diamond(q)$, and
   \item[(b)] $(-1)^{n-1}Z_G(q,\bv)>0$ for all for all  $2$-connected graphs
$G=(V,E)$ with $n$ vertices, at least $m$ edges, and $v_e\in  \scrv_m(q)$
for all $e\in E$.
\end{itemize}
\end{conjecture}

\begin{conjecture}\label{conq>1}
Suppose $1<q\leq 32/27$. Then there exists
a strictly increasing sequence of self-dual intervals $\scrv_m(q)$,
$m\geq 2$, such that
\begin{itemize}
   \item[(a)] $\lim\limits_{m\to\infty} \scrv_m(q)=I_\Diamond(q)$, and
   \item[(b)] $(-1)^{n}Z_G(q,\bv)>0$ for all for all  $2$-connected graphs
$G=(V,E)$ with $n$ vertices, at least $m$ edges, and $v_e\in  \scrv_m(q)$
for all $e\in E$.
\end{itemize}
\end{conjecture}

These conjectures are illustrated in Figure~\ref{fig.finalplot}.

For any fixed $m$, intervals $\scrv_m(q)$
satisfying the hypotheses of Theorem~\ref{thm.0<q<1.block.abstract}
can in principle be found by a finite amount of calculation
(i.e., there are finitely many $m$-edge 2-connected graphs to consider),
but the computations seem rather messy for $m \ge 5$.
For instance, for $m=5$ we have not only the 5-cocycle and the 5-cycle,
but also the triangle with two double edges, its dual $K_4 \setminus e$,
the triangle with one triple edge, and its dual $C_4$ with one double edge.
Indeed, for $m \ge 5$ the interval $I_m$ defined in \reff{def.Im},
which arises by considering only the $m$-cocycle and the $m$-cycle,
{\em cannot}\/ satisfy the hypotheses of Theorem~\ref{thm.0<q<1.block.abstract}
for small $q>0$, as it fails to be contained in $I_\Diamond(q)$:
\begin{eqnarray}
   v_\Diamond^+(q)  & = &  -{q \over 2} \,-\, {q^2 \over 16} \,-\, O(q^3)
      \\[2mm]
   - \, {q \over 1 + (1-q)^{1/m}}   & = &
   -{q \over 2} \,-\, {q^2 \over 4m} \,-\, O(q^3)
\end{eqnarray}

The behaviors expected for $0 < q < 1$ and for $1 < q \le 32/27$
also differ in a curious way.
For $0 < q < 1$ we expect that the upper endpoints $v_m^+(q)$
of the intervals $\scrv_m(q)$ will be strictly increasing towards
$v_\Diamond^+(q)$.
For $1 < q \le 32/27$, by contrast,
we already have $v_m^+(q) = v_\Diamond^+(q)$ {\em exactly}\/
for $m=2,3$ when $9/8 \le q \le 32/27$;
for larger $m$ we can expect this ``crossover point'' $q=9/8$
to move downwards towards $q=1$.
That, at any rate, is our naive guess based on the behavior for small $m$.

Let us note that Conjectures~\ref{conq<1} and \ref{conq>1}
are in a certain sense the most one can hope for,
because we have shown in Corollary~\ref{cor.diamond2a}
that conclusion (b) of Conjecture \ref{conq<1} can only hold
if $q<1$ and $\scrv_m(q)\subseteq I_\Diamond(q)$,
and that conclusion (b) of Conjecture \ref{conq>1} can only hold
if $1<q\leq 32/27$ and $\scrv_m(q)\subseteq I_\Diamond(q)$.

But we can pose the question more broadly,
by asking about the regions
(beyond those covered by Conjectures~\ref{conq<1} and \ref{conq>1})
where we have {\em not}\/ succeeded in controlling the sign of $Z_G(q,\bv)$,
namely:
\begin{itemize}
   \item[(a)] $q < 0$ and $v < -2$;
   \item[(b)] $q < 0$ and $0 < v < -q/2$;
   \item[(c)] $0< q \le 32/27$ ($q \neq 1$) and $v < v_\Diamond^-(q)$;
   \item[(d)] $0< q \le 32/27$ ($q \neq 1$) and $v_\Diamond^+(q) < v < 0$;
   \item[(e)] $q > 32/27$ and $v < 0$.
\end{itemize}
(These regions are labelled ??? in Figure~\ref{fig.finalplot}.)
We {\em conjecture}\/ that if $\scrv$ contains any points in these regions,
then there is {\em no hope}\/ of controlling the sign of $Z_G(q,\bv)$,
at least in terms of the numbers of vertices and edges,
because both signs are possible,
even for the bivariate Tutte polynomial $Z_G(q,v)$:

\begin{conjecture}
   \label{conj.nosigns}
Fix $q$ and $v$ real, and suppose that we are in one of the cases
(a)--(e) above.
Then, for all sufficiently large $n$ (how large depends on $q$ and $v$)
and all sufficiently large $m$ (how large depends on $q$, $v$ and $n$),
there exist 2-connected graphs $G$ with $n$ vertices
and $m$ edges that make $Z_G(q,v)$ nonzero with either sign.
\end{conjecture}

We furthermore suspect that for $0 < q \le 32/27$
(and perhaps all the way up to $q=2$)
the graphs in Conjecture~\ref{conj.nosigns} can be taken to be
series-parallel.
On the other hand, one {\em cannot}\/ use series-parallel graphs
when $q > 2$ and $v \ge -1$, for it is known that $Z_G(q,\bv) > 0$
in this region
\cite[Proposition~6.3 and Corollary~6.4]{Sokal_chromatic_roots}.
More generally, it is known that for graphs of tree-width $\le k$,
we have $Z_G(q,\bv) > 0$ whenever $q > k$ and $v \ge -1$
\cite[Theorem~3.4]{Thomassen_97}.

If Conjecture~\ref{conj.nosigns} is correct,
it follows that the results in this paper,
together with Conjectures~\ref{conq<1} and \ref{conq>1},
are in a fairly strong sense best possible.

We actually conjecture that our results in this paper,
together with Conjectures~\ref{conq<1} and \ref{conq>1},
are best possible in a much stronger sense than that given by
Conjecture~\ref{conj.nosigns}, namely:

\begin{conjecture}
 \label{conj.zeros}
The zeros of the bivariate Tutte polynomials $Z_G(q,v)$
are dense in regions (a)--(e) as $G$ ranges over all graphs.
Here ``dense'' means either that the zeros are dense in $v$
for each fixed $q$, or dense in $q$ for each fixed $v$.
\end{conjecture}

It is possible that Conjectures~\ref{conj.nosigns} and \ref{conj.zeros}
can be proven, at least in cases (c), (d) and (e),
by a variant of the constructions used by
Thomassen \cite[Proposition~2.3 and Theorem~2.5]{Thomassen_97}.

\section*{Acknowledgments}

We wish to thank James Oxley for correspondence concerning
Proposition~\ref{prop.Oxley},
and Tibor Jord\'an for conversations concerning
   Proposition~\ref{prop.combinatorial.graph}
   and its $k$-connected generalization.

We also wish to thank the Isaac Newton Institute for Mathematical Sciences,
University of Cambridge, for generous support during the programme on
Combinatorics and Statistical Mechanics (January--June 2008),
where this work was completed.

This research was supported in part by
U.S.\ National Science Foundation grants PHY--0099393 and PHY--0424082.


\end{document}